\documentclass[reqno]{amsart}

\usepackage{fullpage}
\usepackage{amscd,amssymb}
\usepackage[all]{xy}

\numberwithin{equation}{section}
\newtheorem{Theorem}[equation]{Theorem}

\newtheorem{Lemma}[equation]{Lemma}
\newtheorem{Proposition}[equation]{Proposition}
\newtheorem{Corollary}[equation]{Corollary}

\theoremstyle{definition}

\newtheorem{Definition}[equation]{Definition}
\newtheorem{Remark}[equation]{Remark}
\newtheorem{Example}[equation]{Example}

\newcommand{\BU}[1]{BU\langle #1 \rangle}
\newcommand{\busix}{BU\langle 6 \rangle}
\newcommand{\C}{\mathbb{C}}
\newcommand{\CEll}{\mathcal{C}}
\newcommand{\cp}{\mathbb{C}P^{\infty}}
\newcommand{\cpL}{(\cp)^{L}}
\newcommand{\cpimi}{\cp_{-\infty}}
\DeclareMathOperator*{\colim}{colim}
\DeclareMathOperator{\hocolim}{hocolim}

\newcommand{\e}{0}
\newcommand{\eqdef}{\overset{\text{def}}{=}}
\newcommand{\einfty}{E_{\infty}}
\newcommand{\elly}{Ell_{y}}

\newcommand{\fgpof}[1]{\widehat{#1}}

\newcommand{\Gah}{\widehat{\mathbb{G}}_{a}}
\newcommand{\Gmh}{\widehat{\mathbb{G}}_{m}}
\newcommand{\Gm}{\mathbb{G}_{m}}
\newcommand{\Goc}{G_{\circ}}
\newcommand{\GpOf}[1]{G_{#1}}
\newcommand{\GU}{U_{1}}
\newcommand{\gsharp}{\flat}

\newcommand{\h}{\mathfrak{h}}
\newcommand{\iso}{\cong}
\newcommand{\I}{\mathcal{I}}
\newcommand{\Ioc}{\I_{\circ}}

\newcommand{\JC}{\mathcal{J}_{C}}
\newcommand{\JTate}{\Tate{\mathcal{J}}}
\newcommand{\JEll}{\mathcal{J}}

\DeclareMathOperator{\Ker}{Ker}

\newcommand{\KTate}{\Tate{K}}

\renewcommand{\L}{\mathcal{L}}
\newcommand{\Loo}{\mathcal{A}}
\newcommand{\lsb}[1]{( \! (#1) \! )}
\newcommand{\M}{\mathcal{M}}
\newcommand{\MEll}{M_{\mathrm{Ell}}}
\newcommand{\MU}[1]{MU\langle #1 \rangle}
\newcommand{\musix}{\MU{6}}

\renewcommand{\O}{\mathcal{O}}

\newcommand{\psb}[1]{[ \! [#1] \! ]}
\newcommand{\pihat}{\widehat{\pi}}
\newcommand{\ptit}[1]{#1_{+}}
\newcommand{\ptspace}{*}

\newcommand{\Q}{\mathbb{Q}}

\DeclareMathOperator{\rank}{rank}
\newcommand{\restr}[1]{|_{#1}}
\DeclareMathOperator{\RingSpectra}{RingSpectra}

\newcommand{\slot}{\,-\,}
\newcommand{\Smash}{\wedge}
\DeclareMathOperator{\spec}{spec}

\DeclareMathOperator{\spf}{spf}
\DeclareMathOperator{\sym}      {Sym}

\newcommand{\tensor}[1]{\otimes_{#1}}
\DeclareMathOperator{\Td}{Todd}
\newcommand{\trivialbundle}{\varepsilon}

\newcommand{\Tate}[1]{#1_{\mathrm{Tate}}}
\newcommand{\Tatef}{\Tate{\widehat{C}}}
\newcommand{\TateC}{\Tate{C}}
\newcommand{\tTate}{\Tate{\gamma}}
\newcommand{\T}{\mathbb{T}}

\newcommand{\uln}[1]{\underline{#1}}

\newcommand{\xra}[1]{\xrightarrow{#1}}
\newcommand{\XC}{\mathcal{X}_{C}}

\newcommand{\Z}{\mathbb{Z}}
\newcommand{\ZN}{\T[N]}
\begin{document}

\title[Jacobi orientation]{The Jacobi orientation and the two-variable elliptic genus}

\author[Ando]{Matthew Ando}
\author[French]{Christopher~P.~French}
\author[Ganter]{Nora Ganter}

\address{Department of Mathematics \\ The University of Illinois at
Urbana-Champaign}

\address{Department of Mathematics \\ Grinnell College}

\address{Department of Mathematics \\ The University of Illinois at
Urbana-Champaign}

\email{Ando: mando@math.uiuc.edu}
\email{French: frenchc@math.grinnell.edu}
\email{Ganter: ganter@math.uiuc.edu}

\date{Version 3.1, October 2007}

\begin{abstract}
Let $E$ be an elliptic spectrum with elliptic curve $C$.  We show that
the sigma orientation of \cite{AHS:ESWGTC,ho:icm} gives rise to  
a genus of $SU$-manifolds taking its values in meromorphic functions
on $C$.  As $C$ varies we find that the genus is a meromorphic
arithmetic Jacobi form.  When $C$ is the Tate elliptic curve it
specializes to 
the two-variable elliptic genus studied in 
\cite{EOTY,math-at-0405232,MR91e:57059,DMVV:egspsqs,BL:egsvoeg,BL:egsv,MR2180406}.
We also show that this two-variable genus arises as an instance of the
$S^{1}$-equivariant sigma orientation.
\end{abstract}  

\thanks{Ando was supported by NSF grant DMS-0306429. Ganter was
supported by NSF grant DMS-0504539.  The paper was completed while
Ando and Ganter were visiting MSRI and Stanford University.  We thank
Ralph Cohen for his hospitality.  We thank Alex Ghitza, Johan de Jong,
Michael Hopkins, and Jacob Lurie for useful conversations.}

\maketitle

\section{Introduction}
\label{sec:introduction}

Stable homotopy theory  singles out the Witten genus, also called the
``sigma orientation'', for special attention among elliptic genera.
For example, it is uniquely characterized by
the fact that it refines to a map of $\einfty$ spectra
\[
    MString \rightarrow tmf,
\]
from the bordism spectrum of Spin
manifolds with trivialized second 
Chern class to the spectrum of ``topological modular forms.''  This
spectrum maps canonically and naturally to all elliptic spectra, so
the Witten genus is in this sense ``initial'' among elliptic genera
\cite{ho:icm,AHS:ESWGTC,Hopkins:icm2002}.

On the other hand, the work on orbifold elliptic genera has focused
attention on the two-variable elliptic genus of
\cite{EOTY,MR1189136,math-at-0405232,MR91e:57059}.  This is the genus
for which Dijkgraaf et al. produced a product formula, expressing the
elliptic 
genera of the symmetric product orbifolds $X^{n}/\!/\Sigma_{n}$ as a
function of the 
elliptic genus of $X$ \cite{DMVV:egspsqs}.  Borisov and Libgober have
proved several results about this genus.  For example, they give a
mathematical proof of the product formula.  Most strikingly, they
produce a related ``resolution of singularities'' elliptic genus, and
prove a McKay correspondence result comparing their two 
genera \cite{BL:egsvoeg,BL:egsv,MR2180406}.

In this paper, we explain the relationship between the sigma
orientation and the two-variable elliptic genus.  We express
the relationship two ways.  The first involves the analysis of
$\MU{2p}$-orientations of \cite{AHS:ESWGTC}, and gives new insight on
the modularity properties of the two-variable genus.  The second uses
the sigma orientation in $S^{1}$-equivariant elliptic cohomology, and
gives new insight on the ``level $N$ genera''. 

Let $E$ be a homotopy-commutative even periodic ring spectrum, so $E$
is complex-orientable and 
\[
G=\spf E^{0}\cp
\]
is a (commutative,
one-dimensional) formal group over $\spec \pi_{0}E.$   Let $\MU{2p}$
be the bordism spectrum of manifolds $M$ with complex tangent bundle
and trivializations of $c_{1},\dotsc,c_{p-1},$  so $\MU{2}=MU$ and
$\MU{4}=MSU$.  In
\cite{AHS:ESWGTC}, the authors show that, for $p\leq 3$, the 
set of maps of ring spectra  (or genera, or orientations)
\[
\MU{2p}\to E
\]
is isomorphic to the set of
``$\Theta^{p}$-structures'' on $G$.    (See
\cite{Br:FTTC,AHS:ESWGTC} and \S\ref{sec:thetap-structures})  Briefly,
let $A$ be either an elliptic curve or a one-dimensional formal group,
and let $\I=\I_{A} (\e)$ be the ideal of functions on $A$ vanishing at the
identity.  Let $\Theta^{0}\I=\I,$ and, for $p\geq 1,$ let
$\Theta^{p}\I$ be the line bundle on $A^{p}$ given 
by the formulas 
\begin{align*}
  (\Theta^{1}\I)_{a} & = \frac{\I_{\e}}{\I_{a}} \\
  (\Theta^{2}\I)_{a,b} & =  \frac{\I_{\e} \I_{a+b}}{\I_{a}\I_{b}} \\
 (\Theta^{3}\I)_{a,b,c} & =
\frac{\I_{\e} \I_{a+b}\I_{a+c}\I_{b+c}}{\I_{a}\I_{b}\I_{c}\I_{a+b+c}}, \\
\end{align*}
where $a,b,c$ are points of $A$.  A  $\Theta^{p}$-structure on $A$ is a
trivialization of $\Theta^{p}\I$ which is compatible with various bits
of structure on $\Theta^{p}\I$. 
We write $C^{p} (A;\I)$ for the set of $\Theta^{p}$-structures on $A.$
There are natural maps  
\[
    \delta : C^{p} (A;\I) \rightarrow C^{p+1} (A;\I)
\]
which, in the case of the formal group $\GpOf{E},$ correspond to the
restriction of orientations 
\[
      \delta s: \MU{2p+2}\rightarrow \MU{2p} \xra{s} E.
\]

Notice that there are isomorphisms 
\begin{align*}
   (\Theta^{2}\I)_{a,b}&\iso \Theta^{1}\I_{a}\otimes  \frac{\I_{a+b}}{\I_{b}}\\
   (\Theta^{3}\I)_{a,b,c}&\iso \Theta^{2}\I_{a,b}\otimes 
    \frac{\I_{a+c}\I_{b+c}}{\I_{a+b+c}\I_{c}}.
\end{align*}
As we discuss in \S\ref{sec:sharp}, a section of the right tensor
factor can be viewed as an element of $\O_{X},$ where, in the
$\Theta^{3}$ case for example, $X\subset A^{3}$ is a subspace of 
\[
   \{(a,b,c) \in A^{3} | a+b+c\neq 0,c\neq 0,a+c\neq 0,b+c\neq 0 \}.
\]
For example, if $C$ is an elliptic curve over a ring $R$, then we can take 
\[
   X = \fgpof{C}^{2}\times (C\backslash \e),
\]
by which we mean the pullback of the formal scheme $\fgpof{C}^{2}$ along 
\[
     (C\backslash \e) \to \spec R.
\]
The case of a formal group $G$ over $R$ is more subtle.  In
\S\ref{sec:sharp},  
we consider 
\[
         R_{1} = \colim_{k} \I^{-k} \O_{G};
\]
in terms of a coordinate on $G$ we have
\[
    R_{1} \iso R\lsb{y}.
\]
Let $G_{1}$ and $\Goc$ be formal groups obtained by pulling back $G$
in the diagram  
\[
\begin{CD}
   \Goc @>>>  G_{1} @>>> G \\
    @VVV      @VVV @VVV \\
   \spec R_{1} @>>> \spec \O_{G} @>>> \spec R,
\end{CD}        
\]
(we use $\circ$ to suggest the hole left over when the identity 
is removed from the $G$ in the base).  Then we can use
\[
    X = \Goc^{2}.
\]
To see this, introduce a coordinate on $G$, and let $F$ be the
resulting formal group law.  Then 
\[
   \O_{X} = R\lsb{y}\psb{s,t},
\]
and $y+_{F}s$ and $y+_{F}s+_{F}t$ are units of $\O_{X}.$
We show that a $\Theta^{2}$-structure on $G$ determines a
$\Theta^{1}$-structure on $\Goc,$ and a $\Theta^{3}$-structure
determines a $\Theta^{2}$-structure on $\Goc.$

As we explain in \S\ref{sec:appl-compl-orient}, if $R=\pi_{0}E$ and
$G=\GpOf{E}=\spf 
E^{0}\cp,$ then 
\[
    R_{1} = \pi_{0}E^{\cpimi},
\]
so in  terms of orientations, we find that an $SU$-orientation 
\[
    t: MSU \to E
\]
gives rise to an orientation 
\[
   t^{\sharp}: MU \to E^{\cpimi}.
\]
Starting with a complex orientation $s: MU\to E$, we obtain an
orientation 
\[
 \delta s: MSU \to E
\]
by restriction, and so a new orientation 
\[
   (\delta s)^{\sharp}: MU \to E^{\cpimi}.
\]

This procedure, applied to the Witten genus, produces the two-variable
genus:  let $K\psb{q}$ be the spectrum representing complex $K$-theory with
coefficients extended to $\Z\psb{q},$  and let 
\begin{equation}\label{eq:14}
    \Phi: MU \to K\psb{q}
\end{equation}
be the complex orientation which associates to a manifold $M$ of
complex dimension $d$ the genus
\[
\Td\left(M;\bigotimes_{n\ge 1}\sym_{q^n}(T-\C^{d})
                     \bigotimes_{n\ge 1}\sym_{q^n} (\bar{T}-\C^{d})\right).
\]
Its $K$-theory Euler class is 
\[
   \Phi (u,q)  = (1-u^{-1}) \prod_{n\geq 1}
   \frac{(1-q^{n}u) (1-q^{n}u^{-1})}{(1-q^{n})^{2}};
\]
it is a version of the Witten genus.\footnote{The relationship is
analogous to the relationship between the Todd genus and the $\widehat{A}$
genus.
For example, the genus \eqref{eq:14} coincides with the Witten genus on $SU$-manifolds.} 

As we explain in
\S\ref{sec:sigma-function-two}, the orientation $(\delta \Phi)^{\sharp}$
sends a manifold $M$ of dimension $d$ to the genus 
\begin{equation}\label{eq:26}
\phi (M,y,q) = \Phi (y^{-1},q)^{-d} \Td\left(M;\bigotimes_{n\ge 1}\sym_{q^n}T 
            \bigotimes_{n\ge 1}\sym_{q^n}\bar{T}
            \bigotimes_{n\ge 1}\Lambda_{-yq^n}(T)
            \bigotimes_{n\ge 1}\Lambda_{-y^{-1}q^n}(\bar{T})\right).
\end{equation}
This is one of the standard formulas for the two-variable elliptic
genus, and we give precise comparisons to formulas in
\cite[p. 59]{math-at-0405232} and \cite[p. 4]{BL:egsvoeg}.

Our approach to the two-variable genus gives a new account of its 
modularity, analogous to the account in
\cite{AHS:ESWGTC} of the modularity of the Witten genus.  Abel's
Theorem, or the Theorem of the Cube, implies that an elliptic curve
$C$ has a \emph{unique} $\Theta^{3}$-structure $s (C)$. 
An isomorphism of formal groups
\[
    \gamma: G\iso \hat{C}
\]
then endows $G$ with the $\Theta^{3}$-structure $(\gamma^{3})^{*} s
(C)\restr{\hat{C}^{3}}.$  The data $(E,C,\gamma)$ comprise an
\emph{elliptic spectrum}, and the map of ring spectra 
\[
   s (E,C,\gamma): \musix \to E
\]
arising from the $\Theta^{3}$-structure
is called the \emph{sigma orientation}.

For example, the Tate curve is a (generalized) elliptic curve $\TateC$ over
$\Z\psb{q}$, equipped with an isomorphism  
\[
   \tTate: \Gmh \iso \Tatef.
\]
Now $\Gmh$ is the formal group of complex $K$-theory, so 
$(K\psb{q},\TateC,\tTate)$ is an elliptic spectrum, denoted $\KTate$
for short.   In \cite{Br:FTTC} and \cite[\S2.6,2.7]{AHS:ESWGTC}, it is
shown that\footnote{A \emph{generalized} elliptic curve
may have more than one $\Theta^{3}$-structure, but there is a
unique rule $C\mapsto s (C)$ which is natural in $C$.}
\begin{equation} \label{eq:12}
   s (\TateC) = \delta^{2}\Phi \in \Theta^{3}\I_{\TateC} (\e).
\end{equation}
Now observe that
\[
 s (\TateC)^{\sharp} = (\delta^{2}\Phi)^{\sharp} = \delta
 (\delta\Phi^{\sharp}):     MSU \rightarrow K\psb{q}^{\cpimi},
\]
so the two-variable elliptic genus of an $SU$-manifold is controlled
by the canonical $\Theta^{3}$-structure on the Tate curve.  This leads to a
new proof of the (known) fact that  the two-variable elliptic genus of an
$SU$-manifold is a meromorphic Jacobi form.  To give a precise
statement, let $E\JC$ be the cohomology theory formed by extending the
coefficients of $E$ to $\Gamma (\O_{C\backslash 0}).$  In
\S\ref{sec:arithm-two-vari} we prove the following result. 

\begin{Theorem} 
An elliptic spectrum $(E,C,\gamma)$ determines a \emph{canonical}
orientation of $SU$-manifolds  
\[
   J_{(E,C,\gamma)}: MSU \xra{} E\JC.
\]
The fact that the cubical structure on the Tate curve is of the form
\eqref{eq:12} implies that $J_{\KTate}$ factors through $MU$, and
indeed the diagram 
\[
\begin{CD}
    MSU @> J_{\KTate} >> K\psb{q}\JTate \\
  @VVV @VVV \\
    MU @> \phi >> K\psb{q}\lsb{y}
\end{CD}
\]
commutes, where $\phi$ is the genus \eqref{eq:26}.
\end{Theorem}

We call $J$ the \emph{Jacobi orientation}, because, as we explain in
Theorem \ref{t-th-jacobi}, the naturality of $J$ in the
elliptic spectrum implies that the genus associated to $J_{\KTate}$
takes its values in meromorphic Jacobi forms of index zero.  
We 
emphasize  that $J_{(E,C,\gamma)}$ is canonically 
determined by the elliptic spectrum $(E,C,\gamma)$, and the variety of
expressions for two-variable genera in the literature reflects choices in how to expand it.

Let $\T$ be the circle group.  Jacobi forms also appear as elements of
the $\T$-equivariant elliptic cohomology of spheres of representations.  In
\S\ref{sec:anom-canc-twists}, we give another formula for the
two-variable elliptic genus in terms of the $\T$-equivariant sigma
orientation.    If $V$ is a complex vector bundle over
$X$, let $Vy^{-1}$ denote $V$, considered as $\T$-bundle via the
inverse of the standard action of $\T.$  We also write $y^{-1}$ for
the trivial bundle, with the inverse of the standard action.  Let
$d=\rank T,$ where 
$T$ is the complex tangent bundle of $X$,  and let  
\[
    \xi = Ty^{-1} - T - d y^{-1} = (T-d) (y^{-1}-1)-d.
\]
For $i\geq 0$ let $c_{i}^{\T}$ denote the Borel equivariant Chern
classes.  We note (Lemma \ref{t-le-c-1-c-2-prime-bundle}) that
\[
c_{1}^{\T} (\xi) = 0,
\]
and if $c_{1} (V) = 0$ then 
\[
    c^{\T}_{2} (\xi) = 0.
\]
It follows using \cite{AndoBasterra:WGEEC,Ando:AESO} that if $c_{1}
(V) = 0$, then  $\xi$ has an
equivariant sigma orientation, which is a 
Thom class $U = U (\xi)$ in $E_{\T} (X^{\xi})$, the $\T$-equivariant
elliptic cohomology of Grojnowski associated to a complex elliptic
curve $C$ of the form $\C/ (2\pi i \Z+2\pi i \tau \Z)$ (the first
author and John Greenlees prove an analogous result for Greenlees's
equivariant spectrum \cite{MR2168575} in 
\cite{AG:reso}).  If we compose the Pontrjagin-Thom map
\[
      S^{0}\rightarrow X^{-T}
\]
with the relative zero section 
\[
     X^{-T} \rightarrow  X^{Ty^{-1}-T}
\]
and then desuspend by $dy^{-1}$, we obtain a map 
\begin{equation}\label{eq:24}
   g: S^{-d y^{-1}} \rightarrow X^{\xi}.
\end{equation}
By pulling back the Thom class along this map we obtain an element of 
\[
      \Gamma E_{\T} (S^{-d y^{-1}})\iso \Gamma\I_{C} (\e)^{-d}.
\]
We show (Proposition \ref{t-pr-equiv-const-two-var}) that this is the
two-variable elliptic genus.   We are grateful to Michael Hopkins for
suggesting this approach to the two-variable genus.  The relationship
between our two approaches is explained in Remark \ref{rem-1}.

A similar argument accounts  for the level $N$ elliptic genera of
\cite{MR981372,Witten:Dir,MR1189136}.  If instead of $c_{1} (T) = 0$
we have $c_{1} (T)\equiv 0 \mod N$, then $c_{2}^{\ZN} (\xi) = 0.$  The
work of \cite{AndoBasterra:WGEEC,Ando:AESO} then shows that $\xi$ has
an equivariant sigma orientation $U_{N}\in E_{\ZN} (X^{\xi})$.
Pulling back this Thom class along $g$ as in \eqref{eq:24} gives an
element 
\[
   g^{*}U_{N}\in \Gamma E_{\ZN} (S^{-dy^{-1}})\iso \Gamma \O_{C[N]}.
\]
In Proposition \ref{t-pr-level-N-equiv-const}, we show that the value
of $g^{*}U_{N}$ at $a\in C[N]$ is just the level $N$ elliptic genus of
$X$ evaluated at $a$.

In \cite{MR2254309}, the third author shows that $H_{\infty}$
elliptic genera have a product formula like
that of \cite{DMVV:egspsqs}.  Either of the accounts of the
two-variable elliptic genus given here can be used to prove that the
two-variable genus is $H_{\infty}$, once one knows that the Witten
genus or sigma orientation from which it was constructed is
$H_{\infty}.$  We will return 
to that story at another time.  

We are grateful to the referee for catching a mistake in an earlier
version of the paper and helping us to correct it, and for several
helpful suggestions which improved the paper.  We are 
responsible for the shortcomings which remain.

\section{Notation}
\label{sec:notation}

\subsection{Groups}

We record the notation for some constructions which make sense in any
setting where one has a notion of abelian
group object $G$ over an object $S$, and where the line bundles over
$X$ form a Picard category.  Our examples
will be elliptic curves and formal groups.

If $G$ is an abelian group over $S$, we write $\e: S\to G$ for its identity
section.  If $I$ is a set, we write $G^{I}$ for the product
$G^{I}_{S}$.  If $f: I\to J$, then we write $\pi_{f}$ for the induced
map 
\[
\pi_{f}: G^{J}\to G^{I}.
\]
If $I\subseteq J$ then we may abbreviate this as $\pi_{I}$, and we may
even indicate $I$ by listing its elements.  It will also be convenient
to write $\pihat_{I}$ for $\pi_{J\backslash I}$.  Thus if
$J=\{1,2,3\}$ then $\pi_{13}$ 
and $\pihat_{2}$ indicate the same map $G^{3}\to G^{2}$.  As a special
case we have  
\[
   \pi_{\emptyset} = G^{J} \to S.
\]
We write
\[
\mu_{I}: G^{J} \xra{\pi_{I}} G^{I} \xra{\mu} G
\]
for projection to $G^{I}$ followed by multiplication.  It is
convenient to set 
$\mu_{\emptyset}=\e: G^{J}\to S\to  G$.

\subsection{Change of base for formal groups}

If $G$ is a formal group over $S=\spec A,$
and $T =\spec B\to S$ is an $S$-scheme, then we can form the pull-back
\[
\begin{CD}
G_{T} @>>> G \\
 @VVV @VVV \\
T @>>>  S.
\end{CD}
\]
If $\I (\e)$ denotes the ideal of the identity of $G$, then 
explicitly the ring of functions on $G_{T}$ is the completed tensor product
\[
   \O_{G_{T}} = \left(B \otimes_{A} \O_{G}\right)^{\wedge}_{\I (\e)}.
\]

\subsection{Power series and Laurent series}

Let $R$ be a ring.  Then $R\psb{x}$ will denote the ring of power
series in $x$, and $R\lsb{y}$
will denote the ring of finite-tailed Laurent series in $y.$  In this
paper we will need to consider rings like 
\[
    A = R\lsb{y}\psb{x}.
\]
Note that this is different from  
\[
   B = R\psb{x}\lsb{y}.
\]
For example 
\[
    \sum_{n\geq 0} y^{-n} x^{n} \in A
\]
but 
\[
    \sum_{n\geq 0} y^{-n} x^{n} \not\in B.
\]
An important point is that series of the form 
\[
    y + o (x) \in A
\]
are invertible in $A,$  since $y$ is a unit of $R\lsb{y}$.

Rings such as $R\lsb{y}\psb{x}$ will arise in situations like the
following.  Let $G \iso \spf R\psb{x}$ be a formal group over $\spec
R.$ 
If $G'$ denotes the formal group over $R\lsb{y}$ which is the  pull-back 
\[
\begin{CD}
 G'  @>>> G \\
@VVV @VVV \\
\spec R\lsb{y} @>>> \spec R,
\end{CD}
\]
then 
\[
 \O_{G'} = (R\lsb{y}\otimes_{R} \O_{G})^{\wedge}_{I (\e)} \iso R\lsb{y}\psb{x}.
\]

\subsection{Line bundles}

If $(X,\O_{X})$ is some sort of ringed space, then a \emph{line
bundle} over $X$ will
mean an invertible $\O_{X}$-module, and we write $\L^{\times}$ for the
associated $\Gm$-torsor of trivializations of $\L$;  this participates
in an equivalence of categories between the line bundles and
$\Gm$-torsors over $X$.

\subsection{Vector bundles and Thom spectra}

If $X$ is a space, and $V$ is a vector bundle over $X$, then we write 
$X^{V}$ for the Thom spectrum
\[
X^{V} \eqdef \Sigma^{\infty} (D (V)/S (V)),
\]
where $D (V)$ is the disk
bundle of $V$ and $S (V)$ is the sphere bundle.   
Notice that, if $\trivialbundle$ denotes the trivial bundle of rank
$1$ over $X$, then 
\begin{equation}\label{eq:29}
     X^{\trivialbundle} \iso  \Sigma \Sigma^{\infty} \ptit{X},
\end{equation}
and if $V/X$ 
and $W/Y$, then 
\begin{equation}\label{eq:28}
    (X\times Y)^{V\oplus W} \iso X^{V}\Smash Y^{W}.
\end{equation}

The reason for using spectra rather than spaces is that one can
extend the definition to virtual 
bundles.  For example, if $V$ is a vector bundle over a finite
complex, then there is a vector bundle $W$ over $X$ such that 
\[
     V \oplus W \iso N\trivialbundle.
\]
In view of \eqref{eq:29} and \eqref{eq:28}, it is sensible to set 
\[
    X^{-V} = \Sigma^{-N} X^{W},
\]
and one shows that this stable homotopy type depends only on $V$.
This definition extends to vector bundles over infinite complexes; see
for example \cite{LMS:esht,MR1627486}.

We write $\BU{2p}$ for the connective cover of $\Z\times BU$ with its
bottom non-zero homotopy group in degree $2p.$  So 
\begin{align*}
\BU{0} & = \Z\times BU \\
\BU{2} & = BU \\
\BU{4} & = BSU.
\end{align*}
We write $\MU{2p}$ for the associated Thom spectra, so $\MU{2} = MU$,
and $\MU{0} = MP$ is the two-periodic complex cobordism spectrum.

\subsection{Cohomology}

Let $E$ be a homotopy commutative ring spectrum.  If $X$ is a space,
then $E^{*} (X)$ will denote its \emph{unreduced} cohomology, which is
a ring.   If $Z$ is a spectrum, then $E^{*} (Z)$ will be its usual
spectrum cohomology.  Thus 
\[
     E^{*} (X) = E^{*} (\Sigma^{\infty} \ptit{X}).
\]
With these conventions, a Thom isomorphism in $E$-cohomology for a
vector bundle $V$ of rank $d$ over $X$ is an isomorphism of $E^{*} (X)$-modules
\[
  E^{*} (X) \iso E^{*+d} (X^{V}).
\]
The reduced cohomology of a pointed space $X$ will be written
$\widetilde{E}^{*} (X)$.  

\section{$\Theta^{p}$-structures}
\label{sec:thetap-structures}

We recall from \cite{Br:FTTC,AHS:ESWGTC} the notion of a $\Theta^{p}$-structure
on a line bundle $\L$ over $G$.   For $p\geq 1$ we define line bundles
$\Theta^{p}\L$ over $G^{p}$ by
the formulas
\begin{align}\label{eq:1} 
    \Theta^{p}\L &\eqdef \bigotimes_{I\subseteq \{1,\dotsc ,p \}}
(\mu_{I}^{*}\L)^{(-1)^{|I|}} \\
\end{align}
Thus 
\begin{align*}
  (\Theta^{1}\L)_{a} & = \frac{\L_{\e}}{\L_{a}} \\
  (\Theta^{2}\L)_{a,b} & =  \frac{\L_{\e} \L_{a+b}}{\L_{a}\L_{b}} \\
 (\Theta^{3}\L)_{a,b,c} & =
\frac{\L_{\e} \L_{a+b}\L_{a+c}\L_{b+c}}{\L_{a}\L_{b}\L_{c}\L_{a+b+c}}. \\
\end{align*}
We also define $\Theta^{0}\L=\L$, over $G$.

The formula \eqref{eq:1} for $\Theta^{p}$ exhibits the symmetry of
$\Theta^{p}\L$ over $G^{p}$.    Precisely, we have the following.

\begin{Proposition} \label{t-pr-theta-k-prop}
\begin{enumerate}
 \item For $p>0$, $\Theta^p(\L)$ is a rigid line bundle, that is,
it comes with a trivialization of $\e^{*}\Theta^{p} (\L)$.  
 \item For each permutation $\sigma\in\Sigma_p$, there is a canonical
  isomorphism
  \[ \xi_\sigma:\pi_\sigma^*\Theta^p(\L)\iso\Theta^p(\L). \]
  Moreover, these isomorphisms compose in the obvious way.
 \item There is a canonical identification (of rigid line bundles over
 $G^{p+1}_S$)
 \begin{equation}\label{eq-Theta-cocycle-iso}
  \Theta^p(\L)_{a_1,a_2,\ldots} \otimes
  \Theta^p(\L)_{a_0+a_1,a_2,\ldots}^{-1} \otimes
  \Theta^p(\L)_{a_0,a_1+a_2,\ldots} \otimes
  \Theta^p(\L)_{a_0,a_1,\ldots}^{-1} \iso 1.
 \end{equation}
\end{enumerate}
\end{Proposition}

With these observations, one makes the following

\begin{Definition} \label{defn-cubical-structure}
 Let $\L$ be a line bundle over a group $G$.  A $\Theta^p$--structure on $\L$ is a trivialization $s$ of the line bundle $\Theta^p(\L)$ such that
 \begin{enumerate}
  \item  (rigidity) for $p>0$, $s$ is a rigid section; \label{item-rigid}
  \item  (symmetry) for $p>0$, $s$ is symmetric in the sense that for each
  $\sigma\in\Sigma_p$, we 
   have $\xi_\sigma\pi_\sigma^*s=s$;
  \item (cocycle condition) for $p>1$, the section $s(a_1,a_2,\ldots) \otimes
  s(a_0+a_1,a_2,\ldots)^{-1} 
   \otimes s(a_0,a_1+a_2,\ldots) \otimes s(a_0,a_1,\ldots)^{-1}$ corresponds
   to $1$ under the isomorphism~\eqref{eq-Theta-cocycle-iso}.
 \end{enumerate}
 For simplicity
 a $\Theta^{p}$-structure on $G$ will mean a $\Theta^{p}$-structure on
 the ideal sheaf $\I_{G} (\e).$
 A $\Theta^3$--structure is known as a \emph{cubical structure}
 \cite{Br:FTTC}.  We write $C^p(G;\L)$ for the set of
 $\Theta^p$-structures on $\L$.  Note that $C^0(G;\L)$ is
 just the set of trivializations of $\L$, and $C^1(G;\L)$ is the set
 of rigid trivializations of $\Theta^1(\L)$.   
 Suppose that $G$ begins life as a group over some base $S$.
If $X$  is another object over $S$, then we can write $G_{X}$
 for $G$ considered as a group over the base $X$, and so forth.  In
 that case, we may write $C^{p}(G,X;\L)$ for the set of
 $\Theta^{p}$-structures on  $\L$ over $G_{X}.$
\end{Definition}

Note that $\Theta^{p}$ can also be defined on sections: if $s$ is a
section of $\L$, then $\Theta^{p}s$ is a section of $\Theta^{p}\L$.
In particular if $s$ is a trivialization of $\L$, then $\Theta^{p}s$
is a $\Theta^{p}$-structure on $\L$.

It is important to observe that $\Theta^{p+1}\L$ can be constructed
from $\Theta^{p}\L$ using the group structure in only one factor of
$G^{p}$.  Precisely, if $\M$ is a line bundle over $G\times X$, then we
write $\delta\M$ for 
the line bundle over $G\times G\times X$ given by the formula 
\[
   \delta\M \eqdef \frac{\pi_{1}^{*}\M \pi_{2}^{*}\M}
 {\mu_{12}^{*}\M 
\mu_{\emptyset}^{*}\M}.
\]
That is, 
\[
   \delta \M_{a,b,x} =
      \frac{\M_{a,x}\M_{b,x}}
           {\M_{a+b,x}\M_{\e,x}}.
\]
Let's write $G\vee G$ for the ``wedge''
\[
   G \vee G \eqdef (G\times \{\e \}) \cup (\{\e \}\times G) \subset
   G\times G.
\]
Notice that $(\delta \M)\restr{(G\vee G)\times X}$ is canonically
trivialized, and that from a section $s$ of $\M$ we obtain a section
$\delta s$ of $\delta \M$ in the obvious way.  

\begin{Proposition}\label{t-pr-Theta-delta} 
\begin{enumerate}
\item \label{item:2} For $p\geq 1$, there is a canonical isomorphism of rigid line bundles
\[
   \Theta^{p+1}\L \iso \delta \Theta^{p}\L,
\]
and so 
\[
     \Theta^{p+1}\L \iso \delta^{p}\Theta^{1}\L.
\]
\item \label{item:3} Using this identification, $\delta$ induces a
homomorphism  
\[
     \delta: C^{p} (G^{p},\L) \to C^{p+1} (G^{p+1},\L).
\]
\item \label{item:4} For $p\geq 1$, $\Theta^{p}\L$ is trivialized over the ``fat
wedge'': if 
\[
   i: G^{p-1}\to G^{p}
\]
is any of the $p$ inclusions obtained using the identity of $G$, then
$i^{*}\Theta^{p}\L$ is canonically trivialized.  Moreover,
if $s$ is a $\Theta^{p}$-structure, then $i^{*}s$ coincides with this
trivialization. 
\end{enumerate}
\end{Proposition}

\begin{proof}
Items \ref{item:2} and \ref{item:3} are straightforward.  

For \ref{item:4}, the case $p=1$ is obvious.  For $p>1$, observe that
by symmetry it suffices to treat the case that 
the identity goes to the first 
factor.  In that case
\[
   \Theta^{p} \L_{\e,a_{2}\dotsb} = \delta \Theta^{p-1} \L_{\e,a_{2},\dotsc}
\]
is trivial.  On sections we consider the case $p=2$; the general case
is similar.  Using the cocycle condition we have
\[
     \frac{s (0,0)s (0,b)}{s (0,b)s (0,b)} = 1,
\]
and since $s$ is rigid $s (0,0) = 1$, and so $s (0,b) = 1$ as required.
\end{proof}

\subsection{Complex orientations and $\Theta^{p}$-structures}
\label{sec:compl-orient-thet-1}

We recall how $\Theta^{p}$-structures arise in the
study of multiplicative complex orientations.  The case $p=1$ is the
classical theory of $MU$-orientations, as in 
\cite{Adams:BlueBook}.  The cases $p=2,3$ (and $p=0$) were studied in
\cite{AHS:ESWGTC}.   

An \emph{even periodic} ring spectrum is a ring
spectrum such that $\pi_{\text{odd}}E=0,$ and $\pi_{2}E$ contains a
unit of $\pi_{*}E.$   If $E$ is such a spectrum, then 
\[
   \GpOf{E} \eqdef \spf E^{0}\cp
\]
is a (commutative, one-dimensional) formal group over $S_{E}\eqdef
\spec\pi_{0}E.$  Let $L$ denote the tautological bundle over $\cp$.
The zero section  
\[
    \zeta: \ptit{\cp} \rightarrow (\cp)^{L}
\]
identifies $E^{0} ((\cp)^{L})$ with the (global sections of the) ideal
$\I (\e)$ of functions on $\GpOf{E}$ which vanish at the identity.
The inclusion 
\[
     S^{2} = (\ptspace)^{L}  \rightarrow \cpL
\]
induces isomorphisms
\[
    \pi_{2}E=\widetilde{E}^{0}S^{2} \iso \e^{*}\I (\e) \iso \omega, 
\]
identifying $\pi_{2}E$ with the (global sections) of the sheaf $\omega$
of invariant differentials on $\GpOf{E}.$ 

A map of ring spectra 
\[
    \MU{2p}\to E
\]
gives rise to a map of spectra
\[
    ((\cp)^{p})^{\prod (1-L_{i})} \to \MU{2p} \to E.
\]
The Thom isomorphism in this context can be interpreted as giving a
natural isomorphism of 
\[
E^{0} ((\cp)^{p})=\O_{G_{E}^{p}}
\]
modules
\[
    E^{0} (((\cp)^{p})^{\prod (1-L_{i})}) \iso \Gamma\Theta^{p}\I (\e),
\]
inducing a map 
\[
    \RingSpectra (\MU{2k},E) \to C^{p} (\GpOf{E}^{p};\I (\e)).
\]
About this situation there is the following result of
\cite{AHS:ESWGTC}.

\begin{Theorem} \label{t-th-musix-or}
For $0\leq p\leq 3$ the natural map 
\[
  \RingSpectra (\MU{2p},E) \to C^{p} (\GpOf{E}^{p};\I (\e))
\]
is an isomorphism.  If $s\in C^{p} (\GpOf{E}^{p};\I (\e))$ corresponds
to a map $\MU{2p}\to E,$ then the map 
\begin{equation}\label{eq:13}
    \MU{2p+2}\to \MU{2p} \xra{s} E
\end{equation}
corresponds to $\delta s\in C^{p+1} (\GpOf{E}^{p+1};\I (\e))$.
\qed
\end{Theorem}

\begin{Example}
A map of ring spectra 
\[
   \MU{0} = MP \to E
\]
corresponds to generator $x$ of $\I (\e)$, which is to say an
coordinate on $\GpOf{E}$, or equivalently an element
\[
   U\in E^{0} (\cpL)
\]
whose image $x=\zeta^{*}U$ is a generator of $E^{0}\cp.$
\end{Example}

\begin{Example} \label{ex-3}
A map of ring spectra 
\[
    \MU{2} =  MU\to E
\]
corresponds to a rigid trivialization of $\omega\otimes\I (\e)^{-1}$,
or equivalently of $\omega^{-1}\otimes \I (\e).$  In topology this
corresponds to a dotted arrow making the diagram
\[
\xymatrix{
 {(\cp)^{L-1}}
   \ar@{-->}[r]
&{E}
\\
{S^{0}}
\ar[u]
\ar[ur]
}
\]
commute,
which is the description of complex orientations of $E$ in
\cite{Adams:BlueBook}.  
\end{Example}

Because of Theorem \ref{t-th-musix-or}, if
\[
   s: \MU{2p} \to E
\]
is an orientation, we write $\delta s$ for the induced map \eqref{eq:13}.

\section{Sharp}
\label{sec:sharp}

For $p\geq 1$, a $\Theta^{p+1}$-structure on $G$ nearly defines a
$\Theta^{p}$-structure in the first $p$ variables.  We shall develop
this idea in two ways, but, as we explain in
the introduction, the main point is the following.  For concreteness we
let $G$ be a formal group over a ring $R$, and consider the case
$p=1.$    Let $\I=\I (\e)$ be the ideal sheaf of functions vanishing
at the origin.   In punctual notation, if $a,b$ represent points of
$G$ then   
\[
     (\Theta^{2}\I)_{a,b} \iso
     \frac{\I_{a}}{\I_{\e}}\frac{\I_{a+b}}{\I_{b}} \iso
     (\Theta^{1}\I)_{a} \frac{\I_{a+b}}{\I_{b}}.
\]
In pull-back notation, 
\[
   \Theta^{2}\I \iso \pi_{1}^{*}\left(\Theta^{1}\I \right) \cdot
   \frac{\mu^{*}\I}{\pi_{2}^{*}\I}. 
\]
Now suppose that we have a coordinate on $G$.  Let 
\[
F (x,y) = x + y + O (xy)
\]
be the
resulting formal group law, so 
\[
 \O_{G\times G} \iso R\psb{x,y},
\]
with respect to which 
\[
   \mu^{*}\I = (F (x,y)).
\]

\begin{Lemma} \label{t-le-pull-back-trivial-bc-unit-gen}
$F (x,y)$ is a unit of $R\lsb{y}\psb{x}$,
and so under the ring homomorphism 
\[
    f: R\psb{x,y} \rightarrow S = R\lsb{y}\psb{x},
\]
we have 
\[
    f^{*}\mu^{*}\I  = S.
\]
That is, the ideal sheaf $\mu^{*}\I$ becomes trivial after
pulling back along 
\[
     \spec \left( R\lsb{y} \psb{x} \right) \rightarrow  \spec \O_{G\times G},
\]
and so also over
\[
      \spf \left( R\lsb{y}\psb{x}\right)
\]
Similar remarks hold for the ideal $\pi_{2}^{*}\I.$
\end{Lemma}

\begin{proof}
If we expand $F$ as a power series in $x$, with coefficients power
series in $y$,
\[
   F (x,y) = y + \sum_{i\geq 1} a_{i} (y)  x^{i}
\]
then the constant term $y$ is a unit of $R\lsb{y}.$ 
\end{proof}

In order to take advantage of this observation systematically, 
we introduce the following variant of $\Theta^{p}$.
For $p\geq 0$, let $\Theta^{p}_{*}\L$ be the line bundle over
$G^{p+1}$ given by the formula
\begin{equation}\label{eq:6}
   \Theta^{p}_{*}\L \eqdef \bigotimes_{I\subseteq \{1,\dotsc ,p \}}
(\mu_{I\cup \{p+1\}}^{*}\L)^{(-1)^{|I|+1}}.
\end{equation}


\begin{Example}
In punctual notation, 
\begin{align*}
     (\Theta^{1}_{*}\L)_{a,b} & = \frac{\L_{a+b}}{\L_{b}}\\
     (\Theta^{2}_{*}\L)_{a,b,c} & =
\frac{\L_{a+c}\L_{b+c}}
     {\L_{c} \L_{a+b+c}}.
\end{align*}
\end{Example}

The important relationships between $\Theta^{p}$ and $\Theta^{p}_{*}$
are given by the following result; recall  that 
\[
    \pihat_{p+1}: G^{p+1}\to  G^{p}
\]
denotes projection to the first $p$ factors.

\begin{Proposition}\label{t-pr-theta-st} 
\begin{enumerate}
\item 
For $p>0$, 
\[
      \Theta^{p}\L \iso \left(\Theta^{p}_{*}\L\restr{G^{p}\times\e}
\right)^{-1}. 
\]
\item  (Pascal's Triangle)
\[
      \Theta^{p+1}\L \iso \pihat_{p+1}^{*}\Theta^{p}\L \otimes 
                           \Theta^{p}_{*}\L.
\]
\end{enumerate}
\qed
\end{Proposition}

\begin{Example}
For example 
\begin{align*}
   (\Theta^{2} \L)_{a,b} &= \frac{\L_{0}}{\L_{a}}\frac{\L_{a+b}}{\L_{b}}
\iso   (\Theta^{1}\L)_{a}\otimes (\Theta^{1}_{*}\L)_{a,b} \\
   (\Theta^{3} \L)_{a,b,c} &= \frac{\L_{0}\L_{a+b}}{\L_{a}\L_{b}}
    \frac{\L_{a+c}\L_{b+c}}{\L_{c}\L_{a+b+c}}
\iso   (\Theta^{2}\L)_{a,b}\otimes (\Theta^{2}_{*}\L)_{a,b,c}
\end{align*}
\end{Example}

We are guided by the idea that a section of
$\Theta^{2}_{*}\I (\e)$ restricts to a holomorphic function on the
subspace $X\subseteq G^{3}$ where $a$ and $b$ are small compared to
$c$.   After all, the divisor of
$\Theta^{2}_{*}\I (\e)$ is 
\begin{equation} \label{eq:10}
  [a+c=0] + [b+c = 0] - [c=0] - [a+b+c=0],
\end{equation}
and this divisor intersects $X$ trivially.   In the case of an
elliptic curve $C$ over $S$, we can take 
\[
      X = \fgpof{C}^{2}\times (C\backslash \e).
\]
By this we mean the formal scheme over $U=(C\backslash \e)$ which is the
pull-back of $\fgpof{C}^{2}$ in the diagram  
\[
\begin{CD}
      X @>>> \fgpof{C}^{2} \\
     @VVV      @VVV \\
      U @>>> S;
\end{CD}
\]
in particular $\O_{X}$ is the completed tensor product
\[
  \O_{X} = \O_{\fgpof{C}^{2}} \widehat{\otimes}_{\O_{S}} \O_{U}.
\]

The case of a formal group $G$ over $S=\spec R$ is trickier.  If $y$ is a
coordinate on $G$, then 
\[
   \O_{G} \iso R\psb{y},
\]
with respect to which 
\[
       \I (\e) = (y).
\]
One candidate for $G\backslash \e$ is then ``$\spf R\lsb{y},$'' but 
if $T$ is a discrete $R$-algebra, then the set of \emph{continuous}
maps $R\lsb{y}$ to $T$ 
is the empty set of nilpotent units in $T.$

Let 
\[
   R_{1} = \colim_{k} \I (\e)^{-k} \O_{G} \iso R\lsb{y},
\]
and let $\GU = \spec R_{1}.$  Let $G_{1}$ and $\Goc$ be formal
groups obtained by pulling back $G$ in 
the diagram 
\[
\begin{CD}
\Goc @>>> G_{1} @>>> G \\
@VVV @VVV @VVV \\
\GU @>>> \spec \O_{G} @>>> S.
\end{CD}
\]
By this we mean explicitly that
\[
    \O_{\Goc}\iso R\lsb{y}\psb{x}.
\]
The subscript $1$ on $G_{1}$ indicates lives over a base with one
power series variable.  The subscript $1$ on $\GU$ indicates that
$\O_{\GU}$ has  one Laurent series variable, coming from ``removing
the identity section'' in $\spec \O_{G}.$  The subscript $\circ$
indicates that $G_{\circ}$ lives over a base with a hole in it.

We write $\I$ for the ideal sheaf $\I(\e).$  If it is necessary to
distinguish between $\I_{G} (\e)$ and $\I_{\Goc} (\e)$, then we write
$\Ioc$ for the latter.  

The argument of Lemma \ref{t-le-pull-back-trivial-bc-unit-gen} shows
that $\Theta^{p}_{*}\I$ becomes trivial over $\Goc^{p}$, and so a
section of $\Theta^{p+1}\I\iso \pihat_{p+1}^{*}\Theta^{p}\I \otimes 
                           \Theta^{p}_{*}\I$ gives rise to a section
			   of $\Theta^{p}\Ioc$.  

\begin{Definition}\label{def-sharp} 
If $s$ is a section $\Theta^{p+1}\I$, we write $s^{\sharp}$ for 
resulting section of $\Theta^{p}\Ioc.$ 
\end{Definition}

\begin{Proposition} \label{t-pr-theta-st-over-GU}
If $s$ is a $\Theta^{p+1}$-structure, then $s^{\sharp}$ is a
$\Theta^{p}$-structure, and so 
we have a 
\[
C^{p+1} (G;\I) \to C^{p} (\Goc;\I).
\]
\end{Proposition}
 
\begin{proof}
We need to check the rigidity, symmetry, and cocycle conditions.  
Let $s$ be a $\Theta^{p}$-structure on $\I.$ 

The zero section of $\Goc$ is induced by the map 
\begin{align*}
    G & \xra{} G_{1}^{p} \iso G^{p+1} \\
    c & \mapsto (\e,\dotsc ,\e,c).
\end{align*}
We showed in Proposition \ref{t-pr-Theta-delta} that 
\[
    s (\e,\dotsc ,\e,c) = 1.
\]
The symmetry condition for $s^{\sharp}$ follows easily from the symmetry condition for $s.$

For $p=1$ the cocycle condition is empty. For $p\geq 2$,
the cocycle condition for $s$ does not involve the last variable and 
thus gives the cocycle condition for $s^{\sharp}.$
\end{proof}

\begin{Example}\label{ex-1}
If $s\in C^{2} (G;\I)$, then $s^{\sharp} \in C^{1} (\Goc;\I)$.
Note that then  
\[
     \delta (s^{\sharp}) \in C^{2} (\Goc;\I),
\]
and so we have a homomorphism 
\begin{equation}\label{eq:2}
   C^{2} (G;\I) \to C^{2} (\Goc;\I).
\end{equation}
\end{Example}

\begin{Example}\label{ex-2}
If $t\in C^{3} (G;\I)$, then $t^{\sharp} \in C^{2} (\Goc;\I)$.  In
particular, if $s\in C^{2} (G;\I)$ then  
\[
      (\delta s)^{\sharp} \in C^{2} (\Goc;\I),
\]
and so again we have a homomorphism 
\begin{equation} \label{eq:3}
   C^{2} (G;\I) \to C^{2} (\Goc;\I).
\end{equation}
\end{Example}

\begin{Proposition} \label{t-pr-delta-sharp}
The homomorphisms \eqref{eq:2} and \eqref{eq:3} coincide: for $s\in
C^{2} (G;\I)$ we have 
\[
\delta(s^\sharp )=(\delta s)^\sharp
\]
in $C^{2} (\Goc;\I).$
\end{Proposition}

\begin{proof}
The formula for $\delta s$ involves only the first variable of $s$,
while the construction of $s^{\sharp}$ involves only the last.
\end{proof}

\subsection{$\Theta^{l}$-structures of $\Theta^{k}$-structures}

We describe another approach to the sharp construction which was the
starting point our investigation.  

If $G$ is a formal group over $S$, then we can regard $G^{l}$ as a
group in the first variable, over the base $G^{l-1}.$  In general, we write
$G^{k}_{l-1}$ for the pull-back
\[
\begin{CD}
       G^{k}_{l-1} @>>> G^{k} \\
      @VVV @VVV \\
      G^{l-1} @>>> S.
\end{CD}
\]
Of course we have 
\begin{equation}\label{eq:30}
      G^{k}_{l-1} \iso G^{k+l-1},
\end{equation}
and we shall consider the last $l-1$ factors to be the ``base.''

If $\L$ is a line bundle over $G$, then we can consider $\Theta^{l}\L$
as a line bundle $G_{l-1}.$  We then have two line bundles over
$G^{k}_{l-1},$ namely $\Theta^{k} (\Theta^{l} (\L))$ and
$\Theta^{k+l-1} (\L).$    Explicitly, 
\[ 
\Theta^{k}(\Theta^l(\L)) = \bigotimes_{I\subseteq \{1,\dots,k\}}\mu^*_I(\Theta^l(\L))^{(-1)^{|I|}},
\]
where 
\[
\mu_{I}: G^{k}_{l-1} \iso G^{k+l-1} \to G^{l} \iso G_{l-1}
\]
is the map $\mu_I$ of \S\ref{sec:notation} on the first $k$
factors of $G^{k+l-1}$, and the 
identity on the last $l-1$.
\begin{Proposition}\label{t-pr-kl}
For $k\geq 1$ and $l\geq 2$, there is a canonical isomorphism
\[
\zeta_{l,k} \negmedspace :\Theta^{k+l-1}(\L)^{-1}\cong
\Theta^{k}\left(\Theta^l(\L)\right)
\]
of line bundles over $G^{k+l-1}$.  In particular 
\[
   \Theta^{k+1} (\L)^{-1} \cong \Theta^{k} (\Theta^{2} (\L)).
\]
\end{Proposition}
\begin{proof}
Represent $I\subset \{ 1,\dots ,k\}$ and $J\subset\{ 1,\dots ,l\}$ by vectors 
$(i_1,\dots,i_k)\in\mathbb F_2^k$ and
$(j_1,\dots,j_l)\in\mathbb F_2^l$ respectively.
Then
$$\mu_J\circ\mu_I = \mu_{J\circ I},$$
where $$J\circ I\subset\{1,\dots,k+l-1\}$$ 
denotes the subset represented by 
$$(i_{1}j_{1},\dots ,i_{k}j_{1},j_{2},\dots ,j_{l}) \in\mathbb F_2^{k+l-1}.$$ 
Now
$$\Theta^{k}\Theta^l(\L)=\bigotimes_{I}\mu_I^*\left(\bigotimes_{J}
\mu_J^*(\L)^{(-1)^{|J|+1}}\right)^{(-1)^{|I|}}$$
$$\iso \bigotimes_{I}\bigotimes_{1\in J}\mu_{J\circ
I}^*(\L)^{(-1)^{|I|+|J|+1}},$$ 
because the factors coming from terms with $1\not\in J$ all appear
equally often with  
their inverse and therefore cancel out.
\end{proof}
\begin{Example}\label{ex-ueberfluessig}
We will mainly be interested in the case $k=l=2$.  At a point
$(a,b,c)$ of $G^{2}_{1}\iso G^{3},$ we have 
\begin{align*}
\Theta^{2} (\Theta^{2}\L)_{a,b,c} & \iso 
\frac{(\Theta^{2}\L)_{a+b,c} (\Theta^{2}\L)_{\e,c}}
     {(\Theta^{2}\L)_{a,c} (\Theta^{2}\L)_{b,c}} \\
& \iso \frac{\L_{a+b+c}\L_{a}\L_{b}\L_{c}}{\L_{\e}\L_{a+b}\L_{a+c}\L_{b+c}}\\
& \iso (\Theta^{3}\L)^{-1}.
\end{align*}
Along the same lines, we note that 
\[
  \Theta^{2} (\Theta^{2}\L) \iso \Theta^{2}
  \left(\frac{\mu^{*}\L}{\pi_{1}^{*}\L}\right), 
\]
and in general 
\begin{equation} \label{eq:31}
   \Theta^{p} (\Theta^{2}\L) \iso \Theta^{p}
  \left(\frac{\mu^{*}\L}{\pi_{1}^{*}\L}\right).
\end{equation}
\end{Example}

%
\begin{Definition}
For a $\Theta^{k+1}$-structure $s$
on $\L$, let $s^\gsharp$ be the section
\[
s^\gsharp := \zeta_{2,k}s^{-1}
\]
of $\Theta^k(\Theta^2(\L))$.  
\end{Definition}

Now we can proceed as before.  The proof of Proposition
\ref{t-pr-theta-st-over-GU} applies to give the following. 

\begin{Proposition} \label{t-pr-gsharp}
If $s$ is a 
$\Theta^{k+1}$-structure on $\L$, then $s^\gsharp$ is a
$\Theta^{k}$-structure on 
$\Theta^{2}(\L)$.   \qed
\end{Proposition}
 
A $\Theta^{p+1}$-structure 
\[
s  \in C^{p+1} (G;\I (\e))
\]
gives
rise to a $\Theta^{p}$-structure on $\Theta^{2}\I (\e),$  and indeed
using the isomorphism \eqref{eq:31}, an element
\[
s^{\gsharp} \in C^{p} (G_{1};\mu^{*}\I (\e)/\pi_{1}^{*}\I (\e)).
\]
After pulling back along 
\[
\begin{CD}
     \Goc @>>> G_{1} \\
     @VVV @VVV \\
    \GU @>>> \spec \O_{G},
\end{CD}
\]
we find that the ideal $\mu^{*}\I (\e)$ becomes trivial, and so
$s^{\sharp} = (s^{\gsharp})^{-1}$ can be considered to be an element 
\[
   s^{\sharp}\in C^{p} (\Goc;\I (\e)).
\]

\section{Application to complex orientations}
\label{sec:appl-compl-orient}

In this section we apply the results of \S\ref{sec:sharp} to
orientations of ring spectra.  We begin by describing the topological
counterpart of the base change to $\GU.$

\subsection{The pro-spectrum $\cpimi$} \label{sec:prospectrum-cpimi}

We recall that just as 
\[
    \O_{\GpOf{E}} = E^{0}\cp = \pi_{0}E^{\cp_{+}},
\]
the ring $\O_{\GU}$ arises as the cohomology of Mahowald's pro-spectrum
\[
    \O_{\GU} = E^{0}\cpimi = \pi_{0}E^{\cpimi},
\]
and so $E\O_{\GU} = E^{\cpimi}$.  

In more detail, if $L$ denotes the tautological line bundle over
$\cp$, then $\cp_{-k}$ is the Thom  spectrum of $-k L$, and 
\[
   \cpimi \eqdef (\dotsb \rightarrow \cp_{-k} \rightarrow \cp_{-k+1} \dotsb \rightarrow \cp_{-1} \rightarrow \cp_{+}).
\]
If $t$ is a coordinate on $G$, so $E^{0}\cp= E^{0}\psb{t}$, then 
\[
      E^{0}\cp_{-k} = t^{-k}E^{0}\cp,
\]
and if we define 
\[
      E^{\cpimi} \eqdef  \hocolim  (E^{\cp_{+}} \rightarrow 
                           E^{\cp_{-1}} \rightarrow \dotsb),
\]
then 
\begin{equation} \label{eq:32}
      \pi_{0}E^{\cpimi} \iso E^{0}\lsb{t}.
\end{equation}
It turns out  that
$E^{\cpimi}$ 
is a ring spectrum, in such a way that \eqref{eq:32} is an isomorphism
of rings.   Indeed, the ring structure is represented by a pro-diagonal.  More
precisely, we have the following.

\begin{Proposition}\label{t-pr-cpp} 
There are compatible counit maps
\[
   \cp_{-k} \to S^{0},
\]
and diagonal maps 
\[
\cp_{-k-l} \to \cp_{-k} \wedge \cp_{-l},
\]
such that the obvious
coassociativity, cocommutativity and counit diagrams commute, giving
$\cpimi$ the structure of a comonoid pro-spectrum.
\end{Proposition}

\begin{proof}
The counit map is given by
\[
\cp_{-k} \to \cp_{+} \to S^{0}
\]
which is just the map of Thom spaces associated to the map of bundles 
\[
\begin{CD}
   -kL @>>> \uln{0} \\
   @VVV @VVV \\
   \cp @>>> \ptspace.
\end{CD}
\]
The diagonal is the map of Thom spectra associated to 
\[
\begin{CD}
 -kL - l L @>>> -k L \oplus (-l L)  \\
 @VVV @VVV \\
\cp @>>> \cp \times \cp.
\end{CD}
\]
\end{proof}

Using the ring spectrum structure on $E$, we have compatible maps 
\[
   E^{\cp_{-k}} \Smash E^{\cp_{-l}} \to E^{\cp_{-k}\Smash \cp_{-l}}
   \to E^{\cp_{-k-l}},
\]
and passing to colimits we have the following.  

\begin{Corollary} \label{t-co-cpimi-eou}
If $E$ is an even periodic ring spectrum, then $E^{\cpimi}$ is a
ring spectrum.  
If $G=\spf E^{0}\cp$ is the formal group associated to $E,$ then 
\[
      \pi_{0}E^{\cpimi} \iso \O_{\GU}.
\]
and the formal group associated to $E^{\cpimi}$  is $\Goc.$
A coordinate on
$G$ gives an element $y\in \pi_{0}E^{\cpimi}$ and $x\in E^{0} (\cp)$,
in terms of which 
\[
   \pi_{0}E^{\cpimi}\iso E^{0}\lsb{y},
\]
and 
\[
   (E^{\cpimi})^{0} (\cp) \iso E^{0}\lsb{y} \psb{x}. 
\]
More generally, 
\[
    (E^{\cpimi})^{*} (X) \iso E^{*} (X;E^{*}\lsb{y}).
\] \qed
\end{Corollary}

\begin{Remark}
The fact that $E^{\cpimi}$ is a ring spectrum is well-known.  For
example, this ring structure was studied by Cohen, Jones, and Segal in
\cite{CJS:Floer}.  It can also be deduced from work of Greenlees and
May, for example Proposition 3.5 of \cite{MR1230773}.    Indeed if $E$
is an $\einfty$ ring spectrum, then so is $E^{\cpimi}.$  We shall
address that refinement at another time.
\end{Remark}

\subsection{New orientations from old}
\label{sec:new-orient-from}

Let $E$ be an even periodic ring spectrum with formal group
$G=\GpOf{E}$.  Suppose that we are given a map of 
ring spectra 
\[
    \MU{2p} \to E,
\]
and so a $\Theta^{p}$-structure 
\[
     s \in C^{p} (G;\I (\e)).
\]
Proposition \ref{t-pr-theta-st-over-GU} gives the section 
\[
    s^{\sharp} \in C^{p-1} (\Goc;\I (\e)).
\]
Theorem \ref{t-th-musix-or} and Corollary \ref{t-co-cpimi-eou} imply
that in terms of orientations we have the following result. 

\begin{Proposition}\label{t-pr-q-var-genera}
If $p\leq 4$, then the function
\[
     (\slot)^{\sharp}: C^{p} (G;\I)\to C^{p-1}
     (\Goc;\I) 
\]
determines a function 
\[
     \RingSpectra (\MU{2p},E) \to 
     \RingSpectra (\MU{2 (p-1)},E^{\cpimi}).
\]\qed
\end{Proposition}

\begin{Example}
Thus an $SU$-orientation of a complex-orientable spectrum $E$ gives
rise to a complex orientation of $E^{\cpimi}.$  We thank the
referee for pointing out to us that this natural transformation arises
from a canonical map of ring spectra  
\begin{equation}\label{eq:34} 
MU \to MSU^{\cpimi}.
\end{equation}
Indeed, Adam's theory of complex orientations (described in Example
\ref{ex-3}) says that to give a map of ring spectra \eqref{eq:34} is
equivalent to giving a map 
\[
f : (\cp)^{1-L} \to MSU^{\cpimi}
\]
making the diagram 
\[
\xymatrix{
 {(\cp)^{1-L}}
   \ar@{-->}[r]
&{MSU^{\cpimi}}
\\
{S^{0}}
\ar[u]
\ar[ur]
}
\]
commute.  Our map arises from a compatible family
of maps 
\[
   f_{k}:  (\cp)^{1-L} \Smash \cp_{-k} \iso (\cp\times \cp)^{1 - L - k
   M}  \to MSU,
\]
where we have written $M$ for the tautological bundle over the second
factor $\cp.$  To give such a map, note that $(1-L) (1-M)$ is an
$SU$-bundle, and so we have a map 
\[
    (\cp\times \cp)^{(1-L) (1-M)} \to MSU.
\]
Moreover for $k\geq 1$ the difference
\[
     (1-L) (1 - M) - ( 1 - L - k M) = LM + (k-1)M
\]
is a genuine vector bundle, and so we have the relative zero section
$\zeta$
in the sequence 
\[
      f_{k}: (\cp)^{(1-L)} \Smash \cp_{-k} \iso 
      (\cp\times \cp)^{1 - L - k M}  \xra{\zeta} (\cp \times \cp)^{(1-L)
      (1-M)} \to MSU.
\]
\end{Example}

\begin{Example}
Similarly we can describe the map 
\[
     MSU \to \musix^{\cpimi}
\]
which corresponds to the natural transformation in the Proposition.
Let $V_{k}$ denote the tautological bundle over $BSU (k).$  The bundle
$(k-V_{k}) (1-M)$ over $BSU (k)\times \cp$ is classified by a map 
\[
      BSU (k)\times \cp \to \busix.
\]
The difference
\[
       (k-V_{k}) (1-M) - (k - V_{k} - k M)  = (k-1)M + VM
\]
is a genuine bundle for $k\geq 1$, and so we have a map of Thom
spectra 
\[
      (BSU (k))^{(k-V_{k})} \Smash \cp_{-k} \to (BSU (k)\times
      \cp)^{(k-V_{k}) (1-M)} \to \musix.
\]
Taking adjoints and passing to colimits gives the desired map.
\end{Example}

\subsection{Two-variable genera}
\label{sec:two-variable-genera}

Suppose that $E$ is an even periodic ring spectrum, and let
$G=\GpOf{E}.$  Let 
\[
i:S^{2}\to \cp
\]
denote the inclusion of the bottom cell.  Then 
\[
     E^{0}S^{2} = \pi_{2}E \iso \I/\I^{2} 
\]
is the dual Lie algebra of $G$, and if $f\in E^{0}\cp$ is considered
as a function on $G$, then 
\[
   i^{*} f = df_{\e}.
\]
If $df_{\e}$ is a generator of $\pi_{2}E$, then by Theorem
\ref{t-th-musix-or}, $f$ determines a map
of ring spectra
\[
    MP = \MU{0} \to E,
\]
and so a complex orientation 
\begin{equation}\label{eq:4}
    MU = \MU{2} \to E,
\end{equation}
for which the corresponding element of $C^{1} (G;\I (\e))$ is 
\[
    s = \frac{df_{\e}}{f}.
\]
Again by Theorem \ref{t-th-musix-or}, the $SU$-orientation 
\[
   MSU = \MU{4} \to MU\to  E;
\]
corresponds to the element 
\[
    \delta s = \frac{df_{\e} \mu^{*}f}{\pi_{1}^{*}f\pi_{2}^{*}f}
\]
of $C^{2} (G;\I (\e)).$  

Now Proposition~\ref{t-pr-theta-st-over-GU} implies that 
\[
    (\delta s)^{\sharp} \in C^{1} (\Goc;\I (\e)),
\]
and so determines a complex orientation 
\[
    MU\to E^{\cpimi}.
\]

\begin{Definition}\label{def-two-var} 
Let $E$ be an even-periodic, homotopy-commutative ring spectrum, and let 
\[
\varphi: MU\to E
\]
be a multiplicative complex orientation, associated to a
$\Theta^{1}$-structure $s\in C^{1} (\GpOf{E};\I (\e)).$   The
\emph{adjoint genus} of $\varphi$ 
is the map of ring spectra 
\[
        \widehat{\varphi}: MU \to E^{\cpimi}
\]
associated to the element $(\delta s)^{\sharp} \in C^{1} (\Goc;\I (\e)).$
\end{Definition}

It is illuminating to spell this out in terms of a coordinate $t$ on
$G$, so that $E^{0}\cp\iso E^{0}\psb{t}.$  Then $\O_{G\times G} =
\pi_{0}E\psb{t_{1},t_{2}}$, and the group structure of $G$ can be
expressed as a formal group law 
\[
    \mu^{*} t = F (t_{1}, t_{2}).
\]
We can write 
\[
    f = f (t) \in E^{0}\psb{t}
\]
with $f' (0) \in (\pi_{0}E)^{\times}$, and then $\delta f$ is given by
the expression 
\[
   \delta f (t_{1},t_{2}) = \frac{df_{\e}f (F (t_{1},t_{2}))}{f
   (t_{1})f (t_{2})}. 
\]
The adjoint genus is the genus associated to this expression, 
with $t_{1}$ considered as the coordinate on the group, and
$t_{2}$ considered as an element of $\pi_{0}E^{\cpimi}\iso \pi_{0}E
\lsb{t_{2}}$.

\begin{Example}
If $E$ is rational, then we can choose the coordinate $t$ so that $F$
is the additive group.  Set $x=t_{1}$ and $z=t_{2}$ above, and let
\[
   g (x,z)= \frac{f (x+z)}{f (x)f (z)}.
\]
If $M$ is a complex manifold with total
Chern class 
\[
   c (M) = \prod (1+x_{i})
\]
then by the topological Riemann-Roch theorem (see for example
\cite{MR42:3780,MR1627486}), the adjoint genus of $M$ associated to $f$ is    
\[
   \int_{M} \prod_{i} \frac{x_{i}}{g (x_{i},z)}.
\]
\end{Example}

\section{The sigma function and the two-variable elliptic genus}
\label{sec:sigma-function-two}

We mentioned in the introduction that the Tate elliptic curve $\TateC$
over $\Z\psb{q}$ gives rise to an elliptic spectrum $\KTate$, whose
underlying spectrum is $K\psb{q}.$   In this section we show that when
the analysis in \S\ref{sec:two-variable-genera} is applied to the
complex orientation of $\KTate$ given by the Weierstrass
sigma function, the resulting two-variable genus is the two-variable
elliptic genus. 

\subsection{The sigma function and $\Phi$}\label{sec:sigma-function}
Let $\Phi$ and $\sigma$ be the power series 
\begin{align*}
   \Phi (u,q) & = (1-u^{-1}) \prod_{n\geq 1}
   \frac{(1-q^{n}u) (1-q^{n}u^{-1})}{(1-q^{n})^{2}} \\
   \sigma (u,q) & = u^{1/2}\Phi (u,q) \\
                & = (u^{1/2}-u^{-1/2}) \prod_{n\geq 1}
   \frac{(1-q^{n}u) (1-q^{n}u^{-1})}{(1-q^{n})^{2}}.
\end{align*}
By considering $u$ to be a complex line bundle, one sees
that they define the same genus 
\[
    MSU \rightarrow K\psb{q}.
\]
The genus associated to $\Phi$ factors through $MU$, while the genus
associated to $\sigma$ factors through $MSpin$; as such it is known as
the \emph{Witten genus} \cite{MR1189136}.

We may view $\Phi$ and $\sigma$ as functions of variables $x$ and $\tau$ by
setting 
\begin{align*}
       u^{r} &= e^{r x}\\
       q     &=e^{2\pi i \tau}.
\end{align*}
for $r\in \Q.$  They are variants of the Weierstrass sigma
function.   A number of the following remarks apply to both $\sigma$
and $\Phi$, but for definiteness we focus on $\Phi$.

For $\tau\in \h$, $0<|q|<1$, and so $\Phi$ is a
holomorphic function of $(x,\tau)\in \C\times\h.$  For fixed $\tau$,
$\Phi (x,\tau)$ vanishes to first order when $x$ is a point of the
lattice 
\[
\Lambda = 2\pi i\Z + 2 \pi i \tau \Z,
\]
and has no other
zeroes.  It is not invariant under translation by $\Lambda$ in
$x$; instead we have 
\[
   \Phi (uq^{n}) = (-1)^{n} u^{-n}q^{-n (n+1)/2} \Phi (u).
\]
It follows that $\Phi$ descends to a holomorphic section of the line
bundle
\[
\Loo = \frac{\C^{\times} \times \C}
{(u,v)\sim (uq^{n}, 
v (-1)^{n}u^{n}q^{n (n+1)/2})}
\]
over 
\[
C=\C^{\times}/q^{\Z}\iso \C/\Lambda,
\]
vanishing to first order at the
origin.  As such $\Phi$ is a trivialization of $\Loo\otimes \I
(\e)$, which is to say an isomorphism of line bundles over $C$
\[
     \Loo \iso \I (\e)^{-1}.
\]

\subsection{The adjoint of the Witten genus is the two-variable
elliptic genus} 
\label{sec:two-vari-ellipt-jacobi-cplx}

For the moment let's write $p$ for the projection 
\[
    p: \C\to C.
\]
The classical story of the sigma function implies that if $y$ and $z$
are two points of $\C$, then 
\begin{equation}\label{eq:33}
   W (x,y,z) =  \frac{\Phi (x+y) \Phi (x+z)}
         {\Phi (x) \Phi (x+y+z)},
\end{equation}
considered as a function of $x$, descends to a meromorphic function on
$C$ with divisor 
\[
[-p (y)] + [-p (z)] - [\e ] - [-p (y) - p (z)].
\]

As a function of $x$ we have 
\[
 \Phi (x,\tau) = x + O (x^{2}),
\]
and so via the isomorphism of formal groups 
\[
    \hat{p}: \Gah=\widehat{\C}\to \widehat{C},
\]
$\Phi (x,\tau)$ gives a coordinate on the formal group of $C$.  As
such, 
\[
s  = \frac{\Phi (0,\tau)}{\Phi (x,\tau)}
\]
defines an element of $C^{1} (\fgpof{C},\I (\e))$, and so by Theorem
\ref{t-th-musix-or} determines an orientation 
\[
     MU\to H\Lambda.
\]
Here 
\[
H\Lambda^{*} (X) = HP^{*} (X;\O_{\h}) = H^{*} (X;\O_{\h}[v,v^{-1}])
\]
is periodic ordinary cohomology with coefficients in the holomorphic
function on the upper half plane.    

The resulting $MSU$ orientation corresponds by Theorem
\ref{t-th-musix-or} to the $\Theta^{2}$-structure 
\begin{equation}\label{eq:5}
     \delta s = \frac{\Phi (0,\tau)\Phi (x+z,\tau)}
                     {\Phi (x,\tau)\Phi (z,\tau)}
              = \frac{\sigma (0,\tau)\sigma (x+z,\tau)}
                     {\sigma (x,\tau)\sigma (z,\tau)} \in C^{2}
		     (\fgpof{C}^{2};\I (\e)).  
\end{equation}
As in \S\ref{sec:two-variable-genera}, we then have the adjoint
orientation 
\[
     MU \to H\Lambda^{\cpimi}
\]
associated to the expression \eqref{eq:5}, now written as 
\[
   (\delta s)^{\sharp} \in C^{1} (\fgpof{C}_{\circ};\I (0)).
\]
The associated genus is often called the \emph{two-variable} elliptic genus.

To compare the genus associated to $(\delta s)^{\sharp}$ to standard formulas for the two-variable genus, it is
convenient to use the $q$-expansion formula for $\Phi$ and
express our orientation in $K$-theory, as a map 
\[
   \phi: MU \to K\psb{q}\lsb{y},
\]
where $K\psb{q} (X) = K (X;\Z\psb{q}).$ 
Borisov and Libgober use $-z$ in \eqref{eq:5} where we have used $z$,
and so when passing to $K$-theory it is appropriate to set $y^{r}=e^{-rz}$ for
$r\in \Q$.  We then find that the genus of a manifold $M$ with complex
tangent bundle $T$ of rank $d$ is related to the genus $\elly$ of
\cite[equation (8)]{BL:egsvoeg} by the formula
\begin{equation} \label{eq:16}
\begin{split} 
     \phi (M,y,q) & = \sigma (y^{-1},q)^{-d} 
     y^{-\frac{d}{2}}\Td\left(M;\bigotimes_{n\ge 1} \sym_{q^n}T 
            \otimes \sym_{q^n}\bar{T}
            \otimes \Lambda_{-yq^{n-1}}(\bar{T})
            \otimes \Lambda_{-y^{-1}q^n}(T)\right) \\
                  & =  \sigma (y^{-1},q)^{-d} \elly (M).
\end{split}
\end{equation}
H\"ohn includes the factor of $\sigma (y^{-1},q)$ in the genus.  In
Lemma 2.5.1 of \cite{math-at-0405232} H\"ohn also uses $-z$ where we have
used $z$.  He then sets $y=-e^{z}$, so our $y$ is his $-y$.  With this
understood, we find that 
\begin{equation} \label{eq:15}
\varphi_{\text{H\"ohn}} (M,-y,q) = \phi (M,y,q).
\end{equation}

\subsection{Modularity of the two-variable genus of $SU$-manifolds}

Something interesting happens when we restrict the orientation
$(\delta s)^{\sharp}$ back to $MSU$: this is the orientation 
\[
    MSU \to H\Lambda^{\cpimi}
\]
associated to the section
\begin{equation}\label{eq:9}
  \delta (\delta s)^{\sharp} = \frac{\Phi (0,\tau)\Phi (x+y,\tau)\Phi (x+z)\Phi (y+z)}
                     {\Phi (x,\tau)\Phi (y,\tau) \Phi (z)\Phi
		     (x+y+z)} = (\delta s) W (z,y,x),
\end{equation}
of $\Theta^{2}\I (\e)$ over $\fgpof{C}_{\circ}^{2}$.  We make two
related observations about expression \eqref{eq:9}.
\begin{enumerate}
\item It is precisely the formula for the canonical $\Theta^{3}$-structure
on $\C/\Lambda$, as explained in \cite{Br:FTTC},
\cite[\S2.6]{AHS:ESWGTC}.  
\item Up to the indicated permutation of $x,y,z,$ the factor $W$ in
\eqref{eq:9} is same as the $W$ in \eqref{eq:33}.  Thus
$\delta^{2}s^{\sharp}$ gives in fact a section of
\[
     \Theta^{2}\I (\e)\otimes \mathcal{K}_{\Lambda}, 
\]
where $\mathcal{K}_{\Lambda}$ denotes the meromorphic functions on
$\C/\Lambda.$    As such it determines an orientation 
\[
    MSU \to H\mathcal{K}_{\Lambda}.
\]
\end{enumerate} 

These observations are both complex-analytic aspects of the role of
the Theorem of the Cube in the two-variable genus.  In
\S\ref{sec:arithm-two-vari} we pursue this point of view and construct
a natural genus for $SU$-manifolds taking values in 
meromorphic Jacobi forms; for the curve $\C/\Lambda$ over $\h$
it specializes to give the two-variable genus.   


\section{The Jacobi orientation}
\label{sec:arithm-two-vari}

Abel's Theorem (a particular case of
the Theorem of the Cube) implies that an elliptic curve
$C$ has a canonical cubical structure, that is, an
element $s (C)\in C^{3} (C,\I (\e))$.   If $(E,C,\gamma)$ is an elliptic
spectrum, so $\gamma$ is an isomorphism $\GpOf{E}\iso \fgpof{C},$ then the
``sigma orientation'' of $(E,C,\gamma)$ is the map of ring spectra
\[
\sigma (E,C,\gamma): \MU{6}\to E
\]
associated to
$(\gamma^{3})^{*} (s\restr{\fgpof{C}^{3}})$ by Theorem
\ref{t-th-musix-or} \cite{AHS:ESWGTC}. 

The sigma orientation is \emph{modular} in the following sense.  A
\emph{map} of elliptic spectra 
\[
(f,\alpha): (E,C,\gamma)\to (E',C',\gamma')
\]
is
a map of ring spectra
\[
    f: E\to E'
\]
together with an isomorphism of elliptic curves
\[
      \alpha: C' \iso (\spec \pi_{0}f)^{*}C
\]
making the diagram 
\[
\begin{CD}
  \GpOf{E'}    @> \spf f^{*} >> \GpOf{E} \\
@V \gamma VV @VV \gamma' V \\
  \widehat{C}' @> \widehat{\alpha} >> (\spec \pi_{0}f)^{*}\widehat{C} \\
\end{CD}
\]
commute (recall that $\GpOf{E}= \spf E^{0}\cp$).  Given a map $(f,\alpha)$ of elliptic spectra, the diagram 
\[
\xymatrix{
&
{\musix}
 \ar[dl]_{\sigma_{(E,C,\gamma)}}
 \ar[dr]^{\sigma_{(E',C',\gamma')}}
\\
{E}
 \ar[rr]_{f}
& & 
{E'}
}
\]
commutes.  

As explained in \cite{AHS:ESWGTC}, the preceding discussion extends to
generalized elliptic curves in the sense of \cite{DeligneRapoport}: if
$C/S$ is a generalized elliptic curve, then its smooth locus $C^{reg}$
is a group scheme over $S$, and there is a canonical cubical structure
$s (C)$ on $C^{reg},$   which restricts to a cubical structure on
$\fgpof{C},$ which is a one-dimensional formal group.\footnote{In \cite{AHS:ESWGTC}, a generalized
elliptic curve over $S$ is defined to be a pointed $S$-scheme which is
Zariski locally on $S$ isomorphic to a (possibly singular) 
Weierstrass curve.  The analysis of \cite{Deligne:Tate} discussed
below shows that this definition is equivalent to the definition of
Deligne and Rapoport.}  
The  sigma orientation of the elliptic spectrum associated to the Tate curve is
just the Witten genus \cite{MR1189136,Witten:EllQFT,AHS:ESWGTC}, and the
modularity of the sigma orientation implies that the Witten genus of
an $\musix$-manifold is the $q$-expansion of a modular form.  

In this section we show that the same argument gives 
for any elliptic spectrum $(E,C,\gamma)$ a \emph{canonical, modular} 
$SU$-orientation $J_{(E,C,\gamma)}$, taking values in the spectrum of
meromorphic functions on the curve $C$.
We call it the \emph{Jacobi orientation}.  Its value on the Tate
elliptic spectrum $\KTate$ is the restriction to $MSU$ of the
two-variable elliptic genus $\phi$ 
of \eqref{eq:15}, and its modularity gives a new proof of the fact
that this two-variable genus is a meromorphic weak Jacobi form of index
zero and weight $d$.   

Let $C$ be an elliptic curve over $S=\spec R$, or even a generalized
elliptic curve in the sense of \cite{DeligneRapoport}.   The Jacobi
orientation arises from the  
simple observation that in the isomorphism 
\[
   \Theta^{3}\I_{C} (\e) \iso 
   \Theta^{2}\I_{C} (\e)\otimes \Theta^{2}_{*}\I_{C} (\e)
\]
of line bundles over $C^{3}$, 
the second tensor factor is a trivial subsheaf of the meromorphic functions
on the third factor $C$, by a trivialization which becomes a unit in
$\O_{\fgpof{C}^{2}\times (C\backslash \e)}.$  


Let 
\[
   \JC \eqdef \Gamma (\O_{C\backslash 0}),
\]
and let $\XC = \spec \JC.$  Suppose for simplicity that the formal
group $\fgpof{C}$ admits a coordinate $t$ over
$R.$  Then the
Riemann-Roch Theorem (for a treatment which includes generalized
elliptic curves see \cite{Deligne:Tate}) implies that there are $x \in
\Gamma (\I_{C} (\e)^{-2})$ and $y\in \Gamma (\I_{C} (\e)^{-3})$ and  
$a_{i}\in R$  such that  
\[
    \JC \iso R[x,y]/ (y^{2}+ a_{1}xy + a_{3} y = x^{3} + a_{2} x^{2} +
    a_{4} x + a_{6}) \iso \colim_{k\geq 0} \Gamma (\I_{C} (\e)^{-k}),
\]
and the natural map 
\[
    C\backslash \e \rightarrow  \XC
\]
is an isomorphism.   As in \S\ref{sec:sharp}, let 
\[
   \GU = \spec \colim_{k} \I (0)^{-k} \O_{\fgpof{C}} \iso \spec R\lsb{t}.
\]
Expansion of meromorphic functions at the identity gives a
ring homomorphism 
\[
      \JC \rightarrow \O_{\GU},
\]
and so we have the diagram of formal groups 
\[
\begin{CD}
\fgpof{C}_{\circ} @>>> \fgpof{C}_{\XC} @>>> \fgpof{C}\\
@VVV  @VVV @VVV \\
\GU @>>> \XC @>>> \spec R.
\end{CD}
\]
Then we have the following.


\begin{Proposition}\label{t-pr-cube-to-jacobi}
The canonical cubical structure $s (C) \in C^{3} (C,\I (\e))$
determines a canonical and natural $\Theta^{2}$-structure
\[
    s (C)^{\sharp} \in C^{2} (\fgpof{C},\XC;\I (\e)).
\]
This is compatible with the sharp construction on formal groups of
Proposition \ref{t-pr-theta-st-over-GU} in the sense that 
\[
      s (C)^{\sharp}\restr{\fgpof{C}^{2}_{\circ}} = 
      (s (C)\restr{\fgpof{C}^{3}})^{\sharp} \in 
      C^{2} (\fgpof{C}_{\circ};\I (\e)).
\]
\qed
\end{Proposition}

If $(E,C,\gamma)$ is an elliptic spectrum, then $\fgpof{C}$ admits a
coordinate since 
\[
        \fgpof{C}\iso \spf E^{0}\cp \iso \spf E^{0}\psb{t}.
\]
If we form
the elliptic spectrum $(E\JC,C,\gamma)$, where\footnote{A generalized
elliptic curve $C/R$ is flat over $\spec R$, and so $\JC\iso \Gamma (\O_{C\backslash 0})$
is flat over $R.$} 
\[
    E\JC = E\tensor{\pi_{0}E} \JC,
\]
then we have the $\Theta^{2}$-structure
\[
    (\gamma^{2})^{*} s (C)^{\sharp}
    \in C^{2} (\GpOf{E},\XC;\I (\e)).
\]
Theorem \ref{t-th-musix-or} associates to this $\Theta^{2}$-structure
a multiplicative orientation
\begin{equation} \label{eq:17}
     J_{(E,C,t)}: MSU \rightarrow E\JC.
\end{equation}
Thus we have the following.

\begin{Theorem} 
An elliptic spectrum $(E,C,\gamma)$ determines a canonical map of ring spectra
\[
 J_{(E,C,\gamma)}: MSU\rightarrow E\JC.
\]
Formation of $J_{(E,C,\gamma)}$ is natural, in the sense that if 
\[
    (f,\alpha): (E,C,\gamma) \to (E',C',\gamma')
\]
is a map of elliptic spectra, then the diagram 
\[
\xymatrix{
&
{MSU}
 \ar[dl]_{J_{(E,C,t)}}
 \ar[dr]^{J_{(E',C',t')}}
\\
{E\JC}
 \ar[rr]_{f}
& & 
{E'\mathcal{J}_{C'}}
}
\]
commutes. \qed
\end{Theorem}

\begin{Definition}\label{def-Jacobi} 
The \emph{Jacobi orientation} of the elliptic spectrum $(E,C,\gamma)$ is
the map of ring spectra described by the Theorem. 
\end{Definition}

Let $\MEll$ be the moduli stack of elliptic curves, with universal
curve 
\[
    \pi: \CEll \to \MEll
\]
and identity $\e$.  Let 
\[
\JEll=\mathcal{J}_{\CEll} = \colim_{k\geq 0}  \I (\e)^{-k}
\]
be the indicated sheaf of algebras on $\CEll$, and let 
$\uln{\omega} = \e^{*}\I (\e)$ be the sheaf on $\MEll$  of cotangent
vectors at the origin in $\CEll.$

\begin{Definition}
A \emph{meromorphic weak Jacobi form} of index zero 
and weight $d$ is a global section of 
\[
     \JEll\otimes \pi^{*}\uln{\omega}^{d}.
\]
It is equivalent to give a rule $f$ which associates to each pair
$(C/R,\omega)$, consisting of an elliptic curve $C$ over a ring $R$ and
a trivialization $\omega$ of 
$\uln{\omega}_{C}$ a meromorphic function
\[
    f (C/R,\omega)\in \JC,
\]
subject to the following.
\begin{enumerate}
\item If 
\[
\begin{CD}
C' @> \alpha >> C \\
@VVV @VVV \\
\spec R' @> \alpha >> \spec R
\end{CD}
\]
is a pull-back diagram, then 
\[
      f (C',\alpha^{*}\omega) = \alpha^{*} f (C,\omega).
\]
\item If $\lambda\in R^{\times}$ then 
\[
      f (C,\lambda\omega) = \lambda^{-d} f (C,\omega).
\]
\end{enumerate}
\end{Definition}

\begin{Remark}
This sort of Jacobi form, which might be called an \emph{arithmetic}
Jacobi form, was introduced by Kramer
\cite{MR1310957}.  Its relationship to the usual notion of Jacobi form
as in \cite{MR781735} is the same as the relationship
of the arithmetic to the classical notions of modular form as in 
\cite{Katz:pad}.  
\end{Remark}

\begin{Theorem} \label{t-th-jacobi}
If $X$ is an $SU$-manifold of complex dimension $d$, then as $C$
varies, $J_{(E,C,\gamma)} (X)$ defines a meromorphic Jacobi form of
index zero and 
weight $d$. If $(E,C,\gamma) = \KTate$, and $y$ is the formal function on the
Tate curve corresponding to $u^{-1}$ on $\Gm=\spec\Z[u,u^{-1}],$ then
$J_{\KTate} (X)$ admits an expansion in terms of $y$, and as such
\[
J_{\KTate} (X)_{y} = \phi (X,q,y)
\]
is the restriction to $MSU$ of the two-variable elliptic genus $\phi$
of 
\eqref{eq:16}.  In particular, $\phi (X,q,y)$ is the $q$-expansion of
the meromorphic Jacobi form $J (X).$
\end{Theorem}

\begin{proof}
Let $X$ be an $SU$-manifold of complex dimension $d$.  Then
\[
J_{(E,C,\gamma)} (X) \in \pi_{-2d}E\JC \iso \Gamma
(\uln{\omega}^{d}\otimes \JC),
\]
and so the claim that $J_{(E,C,\gamma)} (X)$ is a Jacobi form follows
from the modularity of the Jacobi orientation, together with the fact
that $\MEll$ has a cover by elliptic spectra.  Alternatively, one can
use the argument in the introduction and \S2.7 of \cite{AHS:ESWGTC} to
show that $\Tate{J} (X)$ is the $q$-expansion of a meromorphic Jacobi form 
of weight $d$ and index $0$.  

The fact that $(J_{\KTate})_{y}$ is the restriction to $MSU$ of the
two-variable elliptic genus follows from the fact that, on the one
hand,  the formula \eqref{eq:9} is the 
$MSU$-characteristic series for the two-variable elliptic genus.  On
the other hand, \eqref{eq:9} is also  the formula for the cubical
structure on a complex elliptic curve of the form $\C/\Lambda$.  Indeed in 
\cite[\S2.6,2.7]{AHS:ESWGTC} this fact is used to show that the
$q$-expansion form of $\sigma$ also gives rise to the cubical
structure on the Tate elliptic curve over $\Z\psb{q}.$
\end{proof}

\section{Anomaly cancellation and twists: the Jacobi genus via 
circle-equivariant elliptic cohomology}
\label{sec:anom-canc-twists} 

In this section, which is independent of
\S\ref{sec:thetap-structures}---\ref{sec:arithm-two-vari},  we show how to obtain the two-variable elliptic genus
of $M$ by calculating the $S^{1}$-equivariant Witten genus of $M$,
twisted by the tangent bundle of $M$, considered as an $S^{1}$-bundle
by the standard action of $S^{1}$ on $TM$.   In fact this method also
leads to an account of the ``level $N$'' genera of 
\cite{MR981372,Witten:Dir,MR1189136}.

\subsection{Umkehr maps and genera}
\label{sec:umkehr-maps-genera}

Let $f: X\to Y$ be a proper map of smooth manifolds.  The
Becker-Gottlieb-Pontrjagin-Thom construction associates to this
situation a stable map  
\[
      \tau (f): \ptit{Y} \to X^{-Tf},
\]
where 
\[
Tf = \Ker df: TX\to f^{*}TY
\]
is the bundle of tangent vectors
along the fiber, and $X^{V}$ is the Thom spectrum of the virtual
vector bundle $V$.   

Associated to any cohomology theory $E$, then, we have a map 
\[
     \tau (f)^{*}: E^{*} X^{-Tf}\to E^{*}Y.
\]
If the bundle $-Tf$ is \emph{oriented} in $E$-theory, that is, we have
a Thom class 
\[
    U \in E^{-d}X^{-Tf},
\]
where $d=\rank Tf$, inducing an isomorphism 
\begin{equation}\label{eq:7}
    E^{*}X \iso E^{*-d}X^{-Tf},
\end{equation}
then the composition of the Thom isomorphism \eqref{eq:7} with $\tau
(f)^{*}$ is the Umkehr homomorphism 
\[
    f_{!}: E^{*}X \to E^{*-d}Y.
\]
We have spelled this out in order to recall the role of the Thom class $U$
in the construction of $f_{!}$, since the notation does not indicate
this dependence.

For example, if $X$ is a compact manifold of dimension $d$, and $T=TX$
is its tangent bundle, then associated to the map 
\[
    \pi^{X}: X\to \ptspace 
\]
is the Pontrjagin-Thom map 
\[
    \tau (\pi): S^{0}\xra{} X^{-T}.
\]
An $E$-orientation of $-T$ gives an Umkehr map 
\[
    \pi_{!}: E^{*}X\to E^{*-d} (\ptspace).
\]
The class $\pi_{!} (1)\in \pi_{d}E$ is the ``genus'' of $X$
associated to the orientation $U$.  Actually the term genus is 
appropriate only in the case that the Thom class $U$ is an instance of
an exponential family of Thom classes, as we now explain.

To give a  map of ring spectra 
\[
    \phi: MU\to E
\]
is equivalent to giving, for every complex vector bundle $V/X$ of rank
$d$, an orientation 
\[
U_{V} \in E^{d} (X^{V}),
\]
which is exponential in the sense that 
\begin{enumerate}
\item if $\epsilon$ is the trivial bundle of rank $1$, then
$U_{\epsilon}$ is the double suspension of $1$ in 
\[
        E^{2} (X^{\epsilon}) = E^{2} (\Sigma^{2}X_{+});
\]
\item if $V/X$ and $W/Y$ are complex
vector bundles of rank $d$ and $e$, then 
\[
   U_{V\oplus W} = U_{V}\Smash U_{W} \in E^{d+e} ((X\times
   Y)^{V\oplus W}) \iso E^{d+e} (X^{V}\Smash Y^{W}).
\]
\end{enumerate}

The effect of $\phi$ on homotopy groups is a ring homomorphism 
\[
   \Phi: MU^{*} (S^{0}) \to E^{*} (S^{0}).
\]
The ring $MU^{*} (S^{0})$ is the bordism ring of manifolds $X$ with
a complex structure on its stable tangent bundle, and so the ring
homomorphism $\Phi$ is a \emph{genus}.  

Let $X$ be such a manifold, of real dimension $2d$.  
Thom's theory \cite{Thom:Cobordism} implies that  
\[
        \Phi (X) = \pi^{X}_{!} (1).
\]

\subsection{Twisted genera}

We can apply $\pi^{X}_{!}$ to classes other than $1$.  For example, if
$\Xi$ is a $E$-characteristic class, then we can define a ``twist'' of 
the genus $\Phi$ as
\[
     \Phi (X;\Xi) \eqdef \pi_{!} (\Xi (TX)).
\]
Or, if $X$ comes equipped with another vector bundle $V$, then we may
form 
\[
    \Phi (X;\Xi (V)) \eqdef \pi_{!} (\Xi (V)).
\]

\subsection{Anomaly cancellation}

A more interesting situation arises when the bundle $-T$ (or
more generally $-Tf$) does \emph{not} admit a Thom class.  
In that case, we may hope to find another vector bundle $V$ on $X$, such that
the virtual bundle $V-Tf$ does admit a
Thom class 
\[
   U \in E^{*} X^{V-Tf}.
\]
If so, then we may compose
the Pontrjagin-Thom map $\tau (f)$ with the 
zero section 
\[
      X^{-Tf} \to X^{V-Tf}
\]
to obtain a stable map 
\[
     \ptit{Y}\to X^{-Tf} \to X^{V-Tf}.
\]
We write 
\[
    f_{!}^{V}: E^{*} (X) \to E^{*} (X^{V-Tf}) \to E^{*} (X^{-Tf}) \to E^{*}Y
\]
for the associated Umkehr map.  In the case of the map $\pi^{X}: X\to
\ptspace$, we obtain a map 
\[
    \pi^{X,V}: E^{*} (X) \rightarrow E^{*}S^{0},
\]
and so we have another kind of twisted genus,
\[
     \Phi (X;V) = \pi^{X,V}_{!} (1).
\]

\subsection{The two-variable elliptic genus as a
specialization of the equivariant Witten genus}

Suppose that $V$ is a complex vector bundle, and let $k$
be an integer.  The reader interested only in the
two-variable genus \eqref{eq:16} can take $k=-1$ in the following.
Let $\T$ be the circle group, let
$y^{k}$ be the one-dimensional complex representation of $\T$ in which
$w\in \T$ acts as $w^{k}$, and let $z=c_{1}^{\T}y\in H^{2}B\T.$  We write
$Vy^{k}$ for $V$, considered 
as a $\T$-equivariant vector bundle using the indicated action of
$\T\subset \C$.  That is, 
\[
    V y^{k} =  V\otimes y^{k}
\]

\begin{Lemma} \label{t-le-c-1-c-2-prime-bundle}
If $V$ is a complex vector bundle, then 
\[
    c_{1}^{\T} (Vy^{k} - V - (\rank V) y^{k}) = 0,
\]
and 
\[
    c_{2}^{\T} (Vy^{k} - V - (\rank V) y^{k}) = -kzc_{1} (V).
\]
In particular, if $c_{1} (V) = 0$, then 
\[
c_{2}^{\T} (Vy^{k} - V - (\rank V) y^{k}) = 0.
\]
\end{Lemma}

\begin{Remark}  \label{rem-1}
The Lemma implies that the map 
\[
   BSU\times \cp \to BSU
\]
classifying $\xi \otimes (L-1)$ factors through $\busix$.  One can see
this using connective $K$-theory, $ku.$  Note that 
\[
   \widetilde{ku}^{2p} (X) \iso [X,\BU{2p}].
\]
So $L-1$ may be viewed as a class 
\[
(L-1) \in  \widetilde{ku}^{2} (\cp),
\]
while the
tautological bundle $\xi$ may be viewed as an element of 
\[
  \xi \in  \widetilde{ku}^{4} (BSU).
\]
Thus 
\[
   \xi\otimes (L-1)\in \widetilde{ku}^{6} (BSU\times \cp) \iso 
[BSU\times \cp,\BU{6}].
\]
The same argument implies that 
\[
    \prod_{i=1}^{p} (1-L_{i}) \in \widetilde{ku}^{2p} ((\cp)^{p}) \iso
    [(\cp)^{p},\BU{2p}], 
\]
which is one of the starting points of \cite{AHS:ESWGTC}.  This sheds
some light on the relationship between our two approaches to the
two-variable elliptic genus.
\end{Remark}

\begin{proof}
Let $d=\rank (V).$  Let 
\[
c = 1 + c_{1} + c_{2}+\dotsb
\]
denote the total
Chern class, and let 
\[
c^{\T} =  1 + c_{1}^{\T} + c_{2}^{\T} + \dotsb 
\]
denote the total Borel Chern class.  If 
\[
  c (V) = \prod (1+x_{i}),
\]
then 
\begin{align*}
  c^{\T} (Vy^{k}) & = \prod (1+ x_{i} + kz)\\
  c^{\T} (d y^{k})  & = (1+kz)^{d} \\
  c^{\T} (Vy^{k} - V - d y^{k}) & = \prod \frac{1+x_{i}+kz}
                                             {(1 + x_{i}) (1+kz)}.
\end{align*}
Without any assumptions about $c_{1} V,$  we have
\[
   c_{1}^{\T} (Vy^{k}) = c_{1}V + d k z = c_{1}^{\T} (V + d y^{k}),
\]
so 
\[
    c_{1}^{\T} (Vy^{k} - V - d y^{k}) = 0.
\]

For $c_{2},$ we find that 
\begin{equation}\label{eq:19}
c_{2}^{\T} (Vy^{k}) = c_{2} V + (d-1) k zc_{1}V + \binom{d}{2}k^{2}z^{2}.
\end{equation}
Taking $V$ to be trivial of rank $d$ in \eqref{eq:19} gives
\[
   c_{2}^{\T} (dy^{k}) = \binom{d}{2}k^{2}z^{2}.
\]
The Whitney sum formula then gives 
\[
c^{\T}_{2} (V + d y^{k}) = c_{2}V + dk z c_{1}V + \binom{d}{2}k^{2}z^{2}.
\]
In general
\begin{equation}\label{eq:8}
   c_{2} (V-W) = c_{2} V - c_{1} V c_{1} W - c_{2}W + c_{1}W^{2},
\end{equation}
which if $c_{1}V=c_{1}W$ simplifies to 
\[
   c_{2} (V-W) =  c_{2}V - c_{2}W.
\]
In our case, this gives 
\begin{align*}
   c_{2}^{\T} (Vy^{k} -  V - dy^{k}) & = c_{2} (V) +
    (d-1)k z c_{1} V    +  \binom{d}{2}k^{2}z^{2} - 
   c_{2} V - d k z c_{1} V - \binom{d}{2}k^{2}z^{2}  \\
& = - k z c_{1}V.
\end{align*}
\end{proof}

We briefly recall some facts about the equivariant elliptic cohomology
theory $E=E_{\T}$ of Grojnowski; for more details see
\cite{Grojnowski:Ell-new,Rosu:Rigidity,AndoBasterra:WGEEC,Ando:AESO}.\footnote{One
can carry out the analysis in this section using Greenlees's 
equivariant elliptic cohomology \cite{MR2168575}; the necessary
prerequisites are the subject of \cite{AG:reso}.}  Let $\Lambda$ be
the lattice  
\[
\Lambda = 2\pi i \Z + 2 \pi i \tau,
\]
and let $C=\C/\Lambda.$   For $a\in C$ let 
\[
    T_{a}: C\to C
\]
denote translation by $a$.  If $X$ is a $\T$-space,
then $E_{\T} (X)$ is a sheaf of $\Z/2$-graded $\O_{\C}$-algebras whose
stalk at a 
point $a\in C$ of exact order $k\leq \infty$ is given by  
\begin{equation}\label{eq:23}
    \left( T_{a}^{*}E_{\T} (X)\right)_{0} = H^{*}_{\T} (X^{\T[k]})\otimes_{H^{*}B\T}\O_{C,0},
\end{equation}
where $\O_{C,0}$ is the stalk of $\O_{C}$ at the identity, and 
\[
z\in
H^{*}B\T\iso \C[z]
\]
is regarded as an element of $\O_{C,0}$ via the projection
\[
   p: \C\to C.
\]
Taking $a=0$ gives 
\[
        E_{\T} (X)^{\wedge}_{0} \iso H\Lambda^{*}_{\T} (X).
\]
We also recall (see for example \cite{Grojnowski:Ell-new} or \cite[Lemma
7.4]{Ando:AESO}) that if $V$ is a complex 
$\T$-vector bundle $V$ over a compact $\T$-space $X$, then 
\[
   E (V)\eqdef E (X^{V})
\]
is an invertible $E (X)$-module.  

The main result of \cite{AndoBasterra:WGEEC} and \cite{Ando:AESO} is
the construction of a Thom class in $E (X^{\xi})$ when $\xi$ is a
virtual $\T$-bundle with $c_{1}^{\T}\xi = 0 = c_{2}^{\T}\xi.$  Applied
to the current situation, their results give the following.

\begin{Proposition}\label{t-pr-equiv-Thom-class}
If
$c_{1}V=0$, then the bundle $Vy^{-1}-V - dy^{-1}$  has a canonical
Thom class   
\[
U\in \Gamma\left( E (Vy^{-1})\otimes E (V)^{-1}\otimes E (y^{-1})^{-d}\right),
\]
whose value in the stalk at the origin is the Borel-equivariant Thom
class given by the sigma orientation of $Vy^{-1}-V-dy^{-1}.$ \qed
\end{Proposition}

In particular, if $V=T=TX$ is the rank-$d$ complex tangent bundle of a 
compact manifold $X$,  then we may consider the composition  
\[
       S^{0} \to X^{-T} \to X^{Ty^{-1}-T},
\]
where the first map is the Pontrjagin-Thom map, and the second is the
relative zero section.  After desuspending by $dy^{-1}$ this gives
\begin{equation} \label{eq:20}
 g: S^{-dy^{-1}}\to X^{-T-dy^{-1}} \to X^{Ty^{-1}-T-dy^{-1}}.
\end{equation}
If $c_{1}X = 0$, then Proposition \ref{t-pr-equiv-Thom-class} gives a
class 
\[
  U \in E (X^{Ty^{-1}-T-dy^{-1}}),
\]
which we may pull back along $g$.  Recall the following.

\begin{Lemma} 
\[
      E (S^{y^{-1}}) = \I (\e),
\]
and so 
\[
      E (S^{-dy^{-1}}) = \I (\e)^{-d}.
\]
\end{Lemma}

\begin{proof}
It's illuminating to give two proofs.  First, consider the cofiber
sequence of $\T$-spaces 
\begin{equation} \label{eq:25}
   \C^{\times} \rightarrow  \C \rightarrow  S^{y^{-1}}.
\end{equation}
Let $S=\spec \C$, and let $\pi$ be the structure map 
\[
   \pi: C\to S.
\]
It's easy to check using \eqref{eq:23} that we have a commutative
diagram 
\[
\begin{CD}
  E_{\T} (\C^{\times}) @<<< E_{\T} (\C) \\
@V \iso VV @VV \iso V \\
  \pi^{*}\O_{S} @< \e^{*} << \O_{C}
\end{CD}  
\]
in which the vertical arrows are isomorphism as indicated.  It follow
that 
\[
   E_{\T} (S^{y^{-1}}) \iso \I (\e),
\]
and the general case follows by taking tensor powers.

Alternatively, observe that \eqref{eq:23} gives 
\[
     \left( T_{a}^{*}E_{\T} (S^{y^{-1}})\right)_{0} \iso \O_{C,0}
\]
for $a\neq \e$, while, letting $L$ denote the line bundle over $\cp$ corresponding to
the representation $y^{-1}$,
\[
    E_{\T} (S^{y^{-1}})_{0} \iso H^{*} ((\cp)^{L})\otimes \O_{C,0}
    \iso \I (\e)_{0}.
\]
This is a description of the line bundle $\I (\e).$
\end{proof}

Let 
\[
\mathcal{E}(X) \eqdef \pi^{T,Ty^{-1}-dy^{-1}}_{!} (1)  = g^{*}U \in E (S^{-dy^{-1}}) =
     \Gamma \I (\e)^{-d}.
\]
\begin{Proposition} \label{t-pr-equiv-const-two-var}
We have 
\[
    \mathcal{E} (X) = J_{(H\Lambda,\C/\Lambda,\hat{p})} (X) = \phi (X,y^{-1},q),
\]
where $\phi$ is the two-variable genus of equation \eqref{eq:16}.
\end{Proposition}

\begin{proof}
Let 
\[
c (T) = \prod_{i=1}^{d} (1+x_{i}).
\]
Let $w$ be a point of $C$, and suppose that $z\in \C$ is such that $p
(z) = w.$ According to \cite[\S6.2]{AndoBasterra:WGEEC} or
\cite[\S8.2]{Ando:AESO}, the class $U$ at $w$ is the Thom class associated to
the Euler class 
\[
    \prod_{i=1}^{d} \frac{\sigma (x_{i}- z,\tau)}{\sigma (x_{i},\tau)\sigma(-z,\tau)},
\]
and we have already shown that this is the  Euler class associated to the
genus $\phi (M,y^{-1},q).$
\end{proof}

\subsection{Level $N$ genera}

If we know only that $c_{1}X\equiv 0 \mod N$, then Lemma
\ref{t-le-c-1-c-2-prime-bundle} shows that 
\[
     c_{2}^{\T[N]} (Ty^{-1} - T - d y^{-1})  = 0.
\]
Noting that 
\[
   E\T\times_{\T} \T\times_{\T[N]} X \iso 
   E\Z/n\times_{\T[N]} X,
\]
we make the following definition.  Again let $E_{\T}$
be Grojnowski's equivariant elliptic cohomology, associated to the
complex curve $C=\C/\Lambda$.

\begin{Definition}\label{def-e-z-n} 
Let $X$ be an $C$-space.  We define the  $\ZN$-equivariant elliptic
cohomology of $X$ to be 
\[
   E_{\ZN} (X) = E_{\T} (\T\times_{\ZN}X),
\]
where $\T\times_{\ZN}X$ is considered as a $\T$ space by acting only
on the left of $\T$.
\end{Definition}

We recall that
\begin{Lemma} \label{t-E-z-n-pt}
\[
E_{\ZN} (\ptspace) = \O_{C[N]}.
\]
More generally,  if $X$ is a $\T[N]$-space, then 
\[
E_{\ZN} (X)_{a} = 0
\]
unless $a\in C[N]$, and if $a$ has exact order $k$ dividing $N$, then 
\[
T_{a}^{*}E_{\ZN} (X)_{a} \iso H^{*}X^{\T[k]},
\]
and the map 
\[
   T_{a}^{*} E_{\T} (\T\times_{\ZN} X)_{a} \iso 
   H_{\T} (\T\times_{\ZN} X^{\T[k]}) \rightarrow 
   H^{*}X^{\T[k]}
\]
corresponds to setting $z=0$ in \eqref{eq:23}.  
\end{Lemma}

\begin{proof}
The stalk of $E_{\T} (Y)$ at a point $a$ of exact order $k\leq \infty$ is 
\begin{equation}\label{eq:18}
    T_{a}^{*}E_{\T} (Y)_{a} = H^{*}_{\T} (Y^{\T[k]})\otimes_{H^{*}B\T} \O_{C,0}.
\end{equation}
If $Y=\T\times_{\ZN} X$ then $Y^{\T[k]}$ is empty unless $k|N$, and
then 
\[
   H^{*}_{\T} (Y^{\T[k]}) = H^{*} (B\ZN\times X^{\T[k]}) = H^{*} (X^{\T[k]})
\]
(recall that we are working with complex coefficients).
\end{proof}

We still have 
\begin{align*}
   c^{\T}_{1} (\T\times_{\ZN} (Ty-T-dy))&  = 0 \\
   c^{\T}_{2} (\T\times_{\ZN} (Ty-T-dy))& = 0
\end{align*}
and so now \cite{AndoBasterra:WGEEC} and \cite{Ando:AESO} imply 

\begin{Proposition}\label{t-pr-level-N-thom-class}
The bundle $Ty^{-1}-T - dy^{-1}$  has a canonical 
Thom class   
\[
   U_{N} \in\left(  E_{\ZN} (Ty^{-1})\otimes E_{\ZN} (T)^{-1}\otimes E_{\ZN} (y^{-1})^{-d} \right).
\]\qed
\end{Proposition}

We can now define 
\[
    \mathcal{E}_{N} (X) = g^{*} U_{N} \in \Gamma E_{\ZN} (S^{-dy^{-1}}).
\]
where, as in \eqref{eq:20}, $g$ is the map 
\[
 g: S^{-dy^{-1}}\to X^{-T-dy^{-1}} \to X^{Ty^{-1}-T-dy^{-1}}.
\]
Now 
\[
   (S^{-y})^{\ZN} = S^{0},
\]
so Lemma \ref{t-E-z-n-pt} implies that 
\[
   E_{\ZN} (S^{-dy^{-1}}) \iso \O_{C[N]}.
\]

\begin{Proposition} \label{t-pr-level-N-equiv-const}
The value of $\mathcal{E}_{N} (X)$ at a point $a\in C[N]$ is the
level-$N$ genus of $X$, as in \cite{MR981372,Witten:Dir},  evaluated at $a$.
\end{Proposition}

\begin{proof}
The recipe for calculating $g^{*}U_{N}$ at $a\in C[N]$ is the
following \cite{Grojnowski:Ell-new}.  Let's write $W$ for our bundle
\[
W = Ty^{-1} - T - d y^{-1}.
\]
Recall that $T_{a}$ is the 
translation map 
\[
   T_{a}: C \to C.
\]
The construction of $E_{\T}$ is such that
\[
 T_{a}^{*}U_{N} \in H^{*}_{\ZN} (M^{W^{A}}).
\]
We may calculate
\begin{equation} \label{eq:22}
    g^{*}T_{a}^{*}U_{N} \in H^{*}_{\ZN} (\ptspace) = \C
\end{equation}
using classical techniques, and this is the value of $\mathcal{E}_{N}$ at $a$.

Let $a = \frac{2\pi i }{N} (l + k \tau),$ with $0\leq k \leq N-1$.
According to \cite[\S9]{Ando:AESO}, $T_{a}^{*}U_{N}$ is the class in
\[
H^{*}_{\T} ((\T\times_{\ZN}X^{W})^{\ZN})\otimes \O_{C,0} \iso 
H^{*}_{\ZN}((\T\times_{\ZN}X^{W})^{\ZN})
\]
whose Euler class in 
\[
   H^{*}_{\T} ((\T\times_{\ZN}X)^{\ZN}) \otimes \O_{C,0}
\]
is
\begin{equation} \label{eq:21}
\exp\left(-\frac{k}{N}\sum x_{i} \right)  
\prod_{i} \frac{\sigma ( x_{i} -z -a)}{\sigma (x_{i}) \sigma (-z-a)}.
\end{equation}
By definition, the quantity $g^{*}T_{a}^{*}U$ in \eqref{eq:22} is the
genus associated to this expression, with $z=0$.  Now observe that
setting $z=0$  in \eqref{eq:21} gives the Euler class associated to
the level-$N$ genus. 
\end{proof}


\begin{thebibliography}{DMVV97}

\bibitem[AB02]{AndoBasterra:WGEEC}
Matthew Ando and Maria Basterra.
\newblock The {W}itten genus and equivariant elliptic cohomology.
\newblock {\em Mathematische Zeitschrift}, 240(4):787--822, 2002,
  arXiv:math.AT/0008192.

\bibitem[Ada74]{Adams:BlueBook}
J.~Frank Adams.
\newblock {\em Stable homotopy and generalised homology}.
\newblock Univ. of Chicago Press, 1974.

\bibitem[AG]{AG:reso}
Matthew Ando and J.~P.~C. Greenlees.
\newblock Circle-equivariant classifying spaces and the rational equivariant
  sigma genus, http://arxiv.org/abs/0705.2687.
\newblock Submitted.

\bibitem[AHS01]{AHS:ESWGTC}
Matthew Ando, Michael~J. Hopkins, and Neil~P. Strickland.
\newblock Elliptic spectra, the {W}itten genus, and the theorem of the cube.
\newblock {\em Inventiones Mathematicae}, 146:595--687, 2001, DOI
  10.1007/s002220100175.

\bibitem[And03]{Ando:AESO}
Matthew Ando.
\newblock The sigma orientation for analytic circle-equivariant elliptic
  cohomology.
\newblock {\em Geometry and Topology}, 7:91--153, 2003, arXiv:math.AT/0201092.

\bibitem[BL02]{BL:egsvoeg}
Lev~A. Borisov and Anatoly Libgober.
\newblock Elliptic genera of singular varieties, orbifold elliptic genus and
  chiral de {R}ham complex.
\newblock In {\em Mirror symmetry, IV (Montreal, QC, 2000)}, volume~33 of {\em
  AMS/IP Stud. Adv. Math.}, pages 325--342. Amer. Math. Soc., Providence, RI,
  2002, arXiv:math.AG/0007126.

\bibitem[BL03]{BL:egsv}
Lev Borisov and Anatoly Libgober.
\newblock Elliptic genera of singular varieties.
\newblock {\em Duke Math. J.}, 116(2):319--351, 2003, arXiv:math.AG/0007108.

\bibitem[BL05]{MR2180406}
Lev Borisov and Anatoly Libgober.
\newblock Mc{K}ay correspondence for elliptic genera.
\newblock {\em Ann. of Math. (2)}, 161(3):1521--1569, 2005,
  arXiv:math.AG/0206241.

\bibitem[Bre83]{Br:FTTC}
Lawrence Breen.
\newblock {\em Fonctions th\^eta et th\'eor\`eme du cube}, volume 980 of {\em
  Lecture Notes in Mathematics}.
\newblock Springer-Verlag, Berlin, 1983.

\bibitem[CJS95]{CJS:Floer}
R.~L. Cohen, J.~D.~S. Jones, and G.~B. Segal.
\newblock Floer's infinite-dimensional {M}orse theory and homotopy theory.
\newblock In {\em The Floer memorial volume}, pages 297--325. Birkh\"auser,
  Basel, 1995.

\bibitem[Del75]{Deligne:Tate}
Pierre Deligne.
\newblock Courbes elliptiques: formulaire (d'apres {J}. {T}ate).
\newblock In {\em Modular functions of one variable {I}{V}}, volume 476 of {\em
  Springer lecture notes}, 1975.

\bibitem[DMVV97]{DMVV:egspsqs}
Robbert Dijkgraaf, Gregory Moore, Erik Verlinde, and Herman Verlinde.
\newblock Elliptic genera of symmetric products and second quantized strings.
\newblock {\em Comm. Math. Phys.}, 185(1):197--209, 1997, arXiv:hep-th/9608096.

\bibitem[DR73]{DeligneRapoport}
P~Deligne and M~Rapoport.
\newblock Les schemas de modules de courbes elliptiques.
\newblock In {\em Modular functions of one variable {I}{I}}, volume 349 of {\em
  Springer lecture notes}, 1973.

\bibitem[Dye69]{MR42:3780}
Eldon Dyer.
\newblock {\em Cohomology theories}.
\newblock Mathematics Lecture Note Series. W. A. Benjamin, Inc., New
  York-Amsterdam, 1969.

\bibitem[EOTY89]{EOTY}
H.~Eguchi, H.~Ooguri, A.~Taormina, and S.-K. Yang.
\newblock Superconformal algebras and string compactification on manifolds with
  ${S}{U}({N})$ holonomy.
\newblock {\em Nucl. Phys. B}, 315, 1989.

\bibitem[EZ85]{MR781735}
Martin Eichler and Don Zagier.
\newblock {\em The theory of {J}acobi forms}, volume~55 of {\em Progress in
  Mathematics}.
\newblock Birkh\"auser Boston Inc., Boston, MA, 1985.

\bibitem[Gan06]{MR2254309}
Nora Ganter.
\newblock Orbifold genera, product formulas and power operations.
\newblock {\em Adv. Math.}, 205(1):84--133, 2006.

\bibitem[GM95]{MR1230773}
J.~P.~C. Greenlees and J.~P. May.
\newblock Generalized {T}ate cohomology.
\newblock {\em Mem. Amer. Math. Soc.}, 113(543):viii+178, 1995.

\bibitem[Gre05]{MR2168575}
J.~P.~C. Greenlees.
\newblock Rational {$S\sp 1$}-equivariant elliptic cohomology.
\newblock {\em Topology}, 44(6):1213--1279, 2005.

\bibitem[Gro07]{Grojnowski:Ell-new}
Ian Grojnowski.
\newblock Delocalized equivariant elliptic cohomology.
\newblock In {\em Elliptic cohomology: geometry, applications, and higher
  chromatic analogues}, volume 342 of {\em London Mathematical Society Lecture
  Notes}. Cambridge University Press, 2007.

\bibitem[H\"oh91]{math-at-0405232}
Gerald H\"ohn.
\newblock Komplexe elliptische {G}eschlechter und ${S}^1$-\"aquivariante
  {K}obordimustheorie (complex elliptic genera and ${S}^1$-equivariant
  cobordism theory), 1991, arXiv:math.AT/0405232.
\newblock Bonn Diplomarbeit.

\bibitem[HBJ92]{MR1189136}
Friedrich Hirzebruch, Thomas Berger, and Rainer Jung.
\newblock {\em Manifolds and modular forms}.
\newblock Aspects of Mathematics, E20. Friedr. Vieweg \& Sohn, Braunschweig,
  1992.
\newblock With appendices by Nils-Peter Skoruppa and by Paul Baum.

\bibitem[Hir88]{MR981372}
Friedrich Hirzebruch.
\newblock Elliptic genera of level {$N$} for complex manifolds.
\newblock In {\em Differential geometrical methods in theoretical physics
  (Como, 1987)}, volume 250 of {\em NATO Adv. Sci. Inst. Ser. C Math. Phys.
  Sci.}, pages 37--63. Kluwer Acad. Publ., Dordrecht, 1988.

\bibitem[Hop95]{ho:icm}
Michael~J. Hopkins.
\newblock Topological modular forms, the {W}itten genus, and the theorem of the
  cube.
\newblock In {\em Proceedings of the International Congress of Mathematicians,
  Vol.\ 1, 2 (Z\"urich, 1994)}, pages 554--565, Basel, 1995. Birkh\"auser.

\bibitem[Hop02]{Hopkins:icm2002}
M.~J. Hopkins.
\newblock Algebraic topology and modular forms.
\newblock In {\em Proceedings of the International Congress of Mathematicians,
  Vol. I (Beijing, 2002)}, pages 291--317, Beijing, 2002. Higher Ed. Press,
  arXiv:math.AT/0212397.

\bibitem[Kat73]{Katz:pad}
Nicholas~M. Katz.
\newblock P-adic properties of modular schemes and modular forms.
\newblock In {\em Modular functions of one variable III}, Lecture Notes in
  Mathematics, pages 70--189. Springer, 1973.

\bibitem[Kra95]{MR1310957}
J.~Kramer.
\newblock An arithmetic theory of {J}acobi forms in higher dimensions.
\newblock {\em J. Reine Angew. Math.}, 458:157--182, 1995.

\bibitem[Kri90]{MR91e:57059}
Igor~M. Krichever.
\newblock Generalized elliptic genera and {B}aker-{A}khiezer functions.
\newblock {\em Mat. Zametki}, 47(2):34--45, 158, 1990.

\bibitem[LMSM86]{LMS:esht}
L.~G. Lewis, Jr., J.~P. May, M.~Steinberger, and J.~E. McClure.
\newblock {\em Equivariant stable homotopy theory}, volume 1213 of {\em Lecture
  Notes in Mathematics}.
\newblock Springer-Verlag, Berlin, 1986.
\newblock With contributions by J. E. McClure.

\bibitem[Ros01]{Rosu:Rigidity}
Ioanid Rosu.
\newblock Equivariant elliptic cohomology and rigidity.
\newblock {\em Amer. J. Math.}, 123(4):647--677, 2001, arXiv:math.AT/9912089.

\bibitem[Rud98]{MR1627486}
Yuli~B. Rudyak.
\newblock {\em On {T}hom spectra, orientability, and cobordism}.
\newblock Springer Monographs in Mathematics. Springer-Verlag, Berlin, 1998.
\newblock With a foreword by Haynes Miller.

\bibitem[Tho54]{Thom:Cobordism}
Renee Thom.
\newblock Quelque proprietes globales des varietes differentiables.
\newblock {\em Comm. Math. Helv.}, 1954.

\bibitem[Wit87]{Witten:EllQFT}
Edward Witten.
\newblock Elliptic genera and quantum field theory.
\newblock {\em Comm. Math. Phys.}, 109, 1987.

\bibitem[Wit88]{Witten:Dir}
Edward Witten.
\newblock The index of the {D}irac operator in loop space.
\newblock In Peter~S. Landweber, editor, {\em Elliptic curves and modular forms
  in algebraic topology}, 1988.

\end{thebibliography}
\def\cprime{$'$}

\end{document}